\documentclass[12pt]{article}
\usepackage{amsmath, amsthm, amsfonts, amssymb,a4wide}
\usepackage{xypic}
\usepackage[all,tips]{xy}

\usepackage{hyperref}
\usepackage[english]{babel}

\usepackage{ifpdf} 
\ifpdf
\usepackage{epstopdf}
\fi

\newtheorem{thm}{Theorem}[section]
\newtheorem{lem}[thm]{Lemma}
\newtheorem{prop}[thm]{Proposition}
\newtheorem{cor}[thm]{Corollary}
\newtheorem{conj}[thm]{Conjecture}

\theoremstyle{definition}
\newtheorem{defin}[thm]{Definition}

\theoremstyle{remark}
\newtheorem{remark}[thm]{Remark}
\newtheorem{example}[thm]{Example}

\newcommand{\bth}{\begin{thm}}
\renewcommand{\eth}{\end{thm}}
\newcommand{\bpr}{\begin{prop}}
\newcommand{\epr}{\end{prop}}
\newcommand{\ble}{\begin{lem}}
\newcommand{\ele}{\end{lem}}
\newcommand{\bco}{\begin{cor}}
\newcommand{\eco}{\end{cor}}
\newcommand{\bde}{\begin{defin}}
\newcommand{\ede}{\end{defin}}
\newcommand{\bex}{\begin{example}}
\newcommand{\eex}{\end{example}}
\newcommand{\bre}{\begin{remark}}
\newcommand{\ere}{\end{remark}}
\newcommand{\bcj}{\begin{conj}}
\newcommand{\ecj}{\end{conj}}

\newcommand{\beq}{\begin{equation}}
\newcommand{\eeq}{\end{equation}}

\newcommand{\ot}{{\otimes}}

\newcommand{\lb}{\label}
\newcommand{\bpf}{\begin{proof}}
\newcommand{\epf}{\end{proof}}
\newcommand{\nl}{\newline}

\newcommand\one           {{\bf1}}

\newcommand{\bD}{{\bold D}}
\newcommand{\bE}{{\bold E}}
\newcommand{\E}{{\cal E}}
\newcommand{\cF}{{\cal F}}

\newcommand{\C}{{\cal C}}

\newcommand{\D}{{\cal D}}

\newcommand{\bC}{{\bold C}}
\newcommand{\I}{{\cal I}}

\newcommand{\Vect}{{\cal V}{\it ect}}

\begin{document}

\vspace*{3em}

\begin{center}
{\LARGE An alternative description of \\[.2em] braided monoidal categories}

\vspace{3em}

{\large 
Alexei Davydov$^{a}$,\, Ingo Runkel$^{b}$,\, ~\footnote{Emails: {\tt alexei1davydov@gmail.com}, {\tt ingo.runkel@uni-hamburg.de}}}
\\[2em]
\it$^a$ 
Department of Mathematics, Ohio University,\\
 Athens, OH 45701, USA
\\[1em]
$^b$ Fachbereich Mathematik, Universit\"at Hamburg\\
Bundesstra\ss e 55, 20146 Hamburg, Germany
\end{center}

\vspace{3em}

\begin{abstract}
We give an alternative presentation of braided monoidal categories. Instead of the usual associativity and braiding we have just one constraint (the b-structure). In the unital case, the coherence conditions for a b-structure are shown to be equivalent to the usual associativity, unit and braiding axioms. We also discuss the next dimensional version, that is, b-structures on bicategories. 
	As an application, we show how special b-categories result in the Yang-Baxter equation, and how special b-bicategories produce Zamolodchikov's tetrahedron equation.
Finally, we define a cohomology theory (the b-cohomology) which plays a role analogous to the one abelian group cohomology has for braided monoidal categories.
\end{abstract}

\newpage

\tableofcontents

\newpage

\section{Introduction}

Consider a not necessarily unital or associative algebra $A$ with multiplication $A\ot A\to A,\ (a,b)\mapsto ab$ such that
\beq \lb{b-algebra-axiom}
	a(bc) = b(ac) ~~ ,\qquad \forall a,b,c\in A \ .
\eeq
We call such algebras {\em b-algebras}. In general, a b-algebra is not even power-associative. 
If, however, a b-algebra is commutative then it is associative:
$$a(bc) = a(cb) = c(ab) = (ab)c \ .$$
This situation occurs for example as a (very) special case of vertex algebras, see \cite[\S\,1.3.3]{fb}.
Next, observe that if a b-algebra is right unital, $a1=a$, then it is commutative:
$$ab = a(b1) = b(a1) = ba \ .$$
Thus right unital b-algebras are nothing but commutative associative unital algebras. In this paper we categorify this observation to arrive at an alternative description of braided monoidal categories.

\medskip

Categorical analogues of associative unital algebras are monoidal categories. 
On the categorical level, commutativity is given by a coherent collection of isomorphisms, called braiding \cite{js}.
Braided monoidal categories have applications in representation theory, low dimensional topology and mathematical physics.
In this paper we categorify b-algebras, and we call the resulting structure a b-category.
As for b-algebras, we find that a right unital b-category is the same as braided monoidal category 
	(Theorem \ref{main}).
Without unit, the b-category structure is much weaker than that of being a braided monoidal category.
Yet, even for non-unital b-categories one can define braid group actions (an important feature of braided monoidal categories). 

The existence of braid group actions is one of our reason to choose the name ``b-category''. Another (related) reason lies in the application to two-dimensional conformal field theory, and specifically in the paper \cite{ms}. There, so-called $B$-matrices are used to express the monodromy of conformal four-point blocks. Similarly, the singular behaviour of four-point blocks as insertion points approach each other is expressed through so-called $F$-matrices. $F$-matrices describe the associator in the category of representations of the corresponding vertex algebra. Our initial motivation was that the monodromy behaviour may be easier to control than the singularities, and so we wanted to understand what categorical structure the $B$-matrices would produce. (However, the results in this paper are purely categorical and no further mention of conformal field theory will be made.)

\medskip

We proceed to categorify once more and obtain b-structures on bicategories. This additional step allows us to make two nice observations in examples we study:

\medskip

\noindent
{\bf\em Yang-Baxter and Zamolodchikov equations:}
Equip the category of vector spaces with the tensor product functor $(V,W) \mapsto V \otimes M \otimes W$ for a fixed vector space $M$. If one analyses b-structures in this case, one finds that they correspond to solutions to the Yang-Baxter equation (Proposition \ref{b-cat-on-vect}). 

If one carries out  a similar construction in a one-object b-bicategory, one is lead to Zamolodchikov's tetrahedron equation \cite{zam} (Theorem \ref{oobb}).
The interpretations of Zamolodchikov's tetrahedron equation as a structure of (a version of) a one-object braided monoidal bicategory was one of the motivations for the study of braided monoidal bicategories \cite{kv}.
We conjecture that a braided monoidal bicategory (in the sense of \cite{bn}) should naturally be a b-bicategory and that this should relate the construction of Zamolodchikov operators on both sides.

Finally, we show that b-functors between such one-object b-bicategories are controlled by an ``$RLLL$-relation'' (Proposition \ref{ibf}) which is also found in three-dimensional integrable lattice models, see e.g.\ \cite{bs}.

\medskip

\noindent
{\bf\em b-cohomology:}
Let us denote a set $S$ with a binary operation satisfying \eqref{b-algebra-axiom} as {\em b-magma}. One can turn a b-magma $S$ into a category by taking its only morphisms to be endomorphisms and defining $S(s,s)$ to be a fixed abelian group $B$ for all $s\in S$. 
If this category is endowed with the structure of a b-category we call it a {\em categorical b-magma}. 
One can carry out an analogous construction to obtain a b-bicategory with objects $S$, only identity 1-morphisms, and a fixed abelian group $B$ for all 2-morphism spaces. 

In this way one is lead to introduce a cohomology theory -- which we call b-cohomology -- whose cochains are maps $S^{\times n} \to B$. One verifies that categorical b-magmas are controlled by $H^3_b(S,B)$ (Proposition \ref{catbmag}) and that the special class of b-bicategories described above is controlled by $H^4_b(S,B)$ (Section \ref{blow}). 

The relation between categorical b-magmas and $H^3_b(S,B)$
is completely analogous to the one between pointed braided monoidal categories (or categorical groups) and abelian (or Eilenberg-Mac Lane) cohomology $H^*_{ab}(A,B)$ for abelian groups \cite{js}. Since every braided monoidal category is a b-category one ends up with a comparison homomorphisms $H^*_{ab}(A,B)$ $\to$ $H^*_b(A,B)$ in degree 2 and 3.  We conjecture that this homomorphism extends to all degrees. 

\bigskip

This paper is organised as follows: In Section \ref{sec:b-categories} we define b-categories and show that unital b-categories are the same as braided monoidal categories. In Section \ref{sec:examples} we give several examples of b-structures. In Section \ref{sec:b-bicategories} we study the categorified version, that is, b-bicategories. In Section \ref{sec:b-cohomology} we introduce b-cohomology and show how it relates to pointed b-categories and b-bicategories.
We conclude in Section \ref{concl} by listing some questions left open by the treatment in this paper.

\subsection*{Acknowledgments}

The authors would like to thank C. Schweigert for suggesting to look at pointed b-categories. AD thanks the Department of Mathematics of Hamburg University for hospitality during the visits in 2011-2013 which were 
supported by the Graduiertenkolleg 1670 of the Deutsche Forschungsgemeinschaft.

\subsection*{Conventions}

We take the definitions of (braided) monoidal categories and functors from \cite{js} as standard with a slight modification - the associativity isomorphism $\alpha$ of a monoidal category and the monoidal structure constraints of a functor, the isomorphisms $\phi_0:I\to F(I)$ and $\phi_2:F(X)\ot F(Y)\to F(X\ot Y)$ from \cite{js}, are replaced by their inverses.

\section{b-categories} \label{sec:b-categories}

\subsection{Definition}

A category $\C$ with a tensor product 
	functor
$$\ot:\C\times \C\to\C,\qquad (X,Y)\mapsto X\ot Y$$
is called {\em b-category} if it comes equipped with a natural in $X,Y,Z\in\C$ collection of isomorphisms
$$\beta_{X,Y,Z}:X\ot(Y\ot Z)\to Y\ot(X\ot Z)\ ,$$
such that the diagram
\beq\lb{cbc}\xymatrix{ & Y\ot(X\ot (Z\ot W)) \ar[rr]^{1\ot\beta_{X,Z,W}} && Y\ot(Z\ot (X\ot W))  \ar[rd]^{\ \beta_{Y,Z,X\ot W}}\\
X\ot(Y\ot (Z\ot W)) \ar[ru]^{\beta_{X,Y,Z\ot W}} \ar[rd]_{1\ot\beta_{Y,Z,W}} &&&& Z\ot(Y\ot (X\ot W))\\
& X\ot(Z\ot (Y\ot W)) \ar[rr]^{\beta_{X,Z,Y\ot W}} && Z\ot(X\ot (Y\ot W)) \ar[ru]_{\ 1\ot\beta_{X,Y,W}}}\eeq
commutes for all $X,Y,Z,W \in \C$.

\medskip

A b-category $\C$ is  {\em symmetric} if $\beta_{X,Y,Z}\beta_{Y,X,Z}=1$ for all $X,Y,Z \in \C$. 

\medskip

A functor $F:\C\to\D$ between b-categories is a {\em b-functor} if it comes equipped with a natural in $X,Y\in\C$ collection of isomorphisms
$$F_{X,Y}:F(X\ot Y)\to F(X)\ot F(Y)\ ,$$
such that the diagram
\beq\lb{cbf}\xymatrix{ 
F(X\ot (Y\ot Z))\ar[rr]^{F_{X,Y\ot Z}} \ar[dd]_{F(\beta_{X,Y,Z})} && F(X)\ot F(Y\ot Z) \ar[rr]^{1\ot F_{Y,Z}} && F(X)\ot (F(Y)\ot F(Z)) \ar[dd]^{\beta_{F(X),F(Y),F(Z)}}\\ \\
F(Y\ot(X\ot Z)) \ar[rr]^{F_{Y,X\ot Z}} && F(Y)\ot F(X\ot Z) \ar[rr]^{1\ot F_{X,Z}} && F(Y)\ot(F(X)\ot F(Z))}\eeq
commutes  for all $X,Y,Z \in \C$.

\medskip

A natural transformation $a:F\to G$ between b-functors $F,G:\C\to\D$ is {\em b-transformation} if the diagram 
\beq\lb{cbt}
\xymatrix{F(X\ot Y) \ar[rr]^{F_{X,Y}}\ar[d]_{a_{X\ot Y}} && F(X)\ot F(Y)\ar[d]^{a_X\ot a_Y}\\ G(X\ot Y)\ar[rr]^{G_{X,Y}} && G(X)\ot G(Y)}\eeq 
commutes for all $X,Y \in \C$.

Composition of b-functors $\xymatrix{\C\ar[r]^G & \D \ar[r]^F & \E}$ is a b-functor with respect to the  isomorphisms $(F\circ G)_{X,Y}:(F\circ G)(X\ot Y)\to (F\circ G)(X)\ot (F\circ G)(Y)$ defined as the composition
\beq\lb{cobf}
\xymatrix{F(G(X\ot Y)) \ar[rr]^{F(G_{X,Y})} && F(G(X)\ot G(Y)) \ar[rr]^{F_{G(X),G(Y)}} && F(G(X))\ot F(G(Y))\ .}
\eeq

\subsection{Braid group action}

Let $X,Y\in\C$ be objects of a b-category. Denote $X^{\ot n}\ot Y=X\ot( X\ot( X\ot (...(X\ot Y)...)$.
Consider the following automorphisms of  
   $X^{\ot n} \ot Y$:
\beq\lb{bgg}b_1=\beta_{X,X,X^{\ot n-2}\ot Y},\qquad b_2=1\ot\beta_{X,X,X^{\ot n-3}\ot Y},\quad ... \quad b_{n-1}=1\ot...\ot 1\ot \beta_{X,X,Y}.\eeq
\bpr
	Let $\C$ be a b-category and $X,Y \in \C$.
Sending the Coxeter generators of the braid group $B_{n}$ to $b_i$ defines a group homomorphism
$$B_{n}\to Aut_\C(X^{\ot n}\ot Y) \ .$$
	If $\C$ is symmetric, this homomorphism factors through the symmetric group $S_n$.
\epr
\bpf
The Coxeter relations $b_ib_{i+1}b_i=b_{i+1}b_ib_{i+1}$ follow from the axiom \eqref{cbc}.
The commutativity relations $b_ib_j=b_jb_i$ for $|i-j|>1$ follow from the naturality of $\beta$.
Clearly, $b_i^2=1$ for symmetric $\C$.
\epf

Similarly one defines the pure braid group action $P_{n}\to Aut_\C(X_1\ot ...\ot X_n\ot Y)$,
where $X_i,Y\in\C$ and $X_1\ot ...\ot X_n\ot Y = X_1\ot( X_2\ot( X_3\ot (...(X_n\ot Y)...)$.

\subsection{Unital b-categories}

A b-category $\C$ is  {\em unital} if there is an object $I\in\C$ (the {\em unit} object) together with a natural in $X\in\C$ collection of isomorphisms
$$\rho_X:X\ot I\to X\ ,$$
	subject to conditions \eqref{ubc3} and \eqref{ubc4} below. To formulate the conditions succinctly, define the collection of isomorphisms $c_{X,Y}:X\ot Y\to Y\ot X$ as
\beq\lb{br}
c_{X,Y} = \Big(
\xymatrix{X\ot Y\ar[rr]^{1\ot\rho^{-1}_Y} && X\ot (Y\ot I) \ar[rr]^{\beta_{X,Y,I}} && Y\ot (X\ot I) \ar[rr]^{1\ot\rho_X}& & Y\ot X}
\Big) \ .
\eeq
The data $I,\rho$ have to make the two diagrams
\beq\lb{ubc3}
\xymatrix{
& (X Y) (Z W) \ar[rr]^{c_{X Y,Z W}} \ar[ld]_{\beta_{X Y,Z,W}} && (Z W) (X Y) \ar[rd]^{\beta_{Z W,X,Y}} \\
Z((X Y) W) \ar[d]^{1 c_{X Y,W}} & && & X((Z W) Y) \ar[d]_{1 c_{Z W,Y}} \\
Z(W(X Y)) \ar[dr]_{1 \beta_{W,X,Y}} & && & X(Y (Z W)) \\
& Z(X(W Y)) \ar[dr]_{\beta_{Z,X,W Y}} & & X(Z (Y W)) \ar[ur]_{1\beta_{Z,Y,W}} \\
&& X(Z(W Y)) \ar[ru]_{1 1 c_{W,Y}} }
\eeq
\beq\lb{ubc4}
\xymatrix{
X(IY) \ar[rr]^{\beta_{X,I,Y}} && I(XY)\\
X(YI) \ar[u]^{1c_{Y,I}} \ar[dr]^{1\rho_Y} && (XY)I \ar[u]_{c_{XY,I}} \ar[dl]_{\rho_{XY}} \\
& XY }
\eeq
commutative for all $X,Y,Z,W \in \C$.

\ble\lb{brc}
The natural collection \eqref{br} in a unital b-category fits into the following commutative diagrams: 
\beq\lb{brc1}\xymatrix{ & Y\ot(X\ot Z) \ar[rr]^{1\ot c_{X,Z}} && Y\ot(Z\ot X)  \ar[rd]^{\ \beta_{Y,Z,X}}\\
X\ot(Y\ot Z) \ar[ru]^{\beta_{X,Y,Z}} \ar[rd]_{1\ot c_{Y,Z}} &&&& Z\ot(Y\ot X)\\
& X\ot(Z\ot Y) \ar[rr]^{\beta_{X,Z,Y}} && Z\ot(X\ot Y) \ar[ru]_{\ 1\ot c_{X,Y}}}\eeq
\beq\lb{ubc2}
\xymatrix{
X\ot (Y\ot Z) \ar[rr]^{c_{X,Y\ot Z}} \ar[dd]_{\beta_{X,Y,Z}} && (Y\ot Z)\ot X \ar[rr]^{c_{Y\ot Z,X}} && X\ot (Y\ot Z)\\ \\
Y\ot(X\ot Z) \ar[rr]^{1\ot c_{X,Z}} && Y\ot(Z\ot X) \ar[rr]^{1\ot c_{Z,X}} && Y\ot(X\ot Z) \ar[uu]_{\beta_{Y,X,Z}}
}
\eeq
\ele
\bpf
Commutativity of \eqref{brc1} is implied by commutativity of the following diagram:
$$\xymatrix{
Y\ot(X\ot Z) \ar[rrrr]^{1\ot c_{X,Z}} &&&& Y\ot(Z\ot X)  \ar[dd]^{\ \beta_{Y,Z,X}}\\
& Y\ot(X\ot (Z\ot I)) \ar[rr]^{1\ot\beta_{X,Z,I}} \ar[ul]_{1(1\rho_Z)} && Y\ot(Z\ot (X\ot I)) \ar[ur]^{1(1\rho_X)} \ar[d]_{\beta_{Y,Z,X\ot I}}\\
X\ot(Y\ot Z) \ar[uu]^{\beta_{X,Y,Z}} \ar[dd]_{1\ot c_{Y,Z}} & X\ot(Y\ot (Z\ot I)) \ar[l]_(.55){1(1\rho_Z)} 
	\ar[d]^{1 \ot \beta_{Y,Z,I}} 
\ar[u]_{\beta_{X,Y,Z\ot I}} && Z\ot(Y\ot (X\ot I)) \ar[r]^(.56){1(1\rho_X)} & Z\ot(Y\ot X)\\
& X\ot(Z\ot (Y\ot I)) 
   \ar[rr]_{\beta_{X,Z,Y \ot I}} 
\ar[dl]_{1(1\rho_Y)} && Z\ot(X\ot (Y\ot I)) \ar[u]^{1\ot \beta_{X,Y,I}} \ar[dr]^{1(1\rho_Y)} \\
X\ot(Z\ot Y) \ar[rrrr]_{\beta_{X,Z,Y}} &&&& Z\ot(X\ot Y) \ar[uu]_{1\ot c_{X,Y}}}$$
Commutativity of \eqref{ubc2} is implied by the following pasting of the coherence \eqref{ubc3} for $X,I,Y,Z$:
$$\xymatrix{
& (X I) (Y Z) \ar[rr]^{c_{X I,Y Z}} \ar[ld]_{\beta_{X I,Y,Z}} \ar[rd]^{\rho_X11} && (Y Z) (X I) \ar[rd]^{\beta_{Y Z,X,I}} \ar[d]_{11\rho_X}\\
Y((X I) Z) \ar[d]_{1 c_{X I,Z}} \ar[r]^{1\rho_X1} & Y(XZ) \ar[d]^{1c_{X,Z}} & X(YZ) 
    \ar[r]^{c_{X,YZ}}  \ar[l]_{\beta_{X,Y,Z}} 
& (YZ)X \ar[d]_{c_{YZ,X}} & X((Y Z) I) \ar[d]^{1 c_{Y Z,I}} \ar[ld]_{1\rho_{YZ}} \\
Y(Z(X I)) \ar[dr]_{1 \beta_{Z,X,I}} \ar[r]^{11\rho_X} & Y(ZX) \ar[rd]^{1c_{Z,X}} && X(YZ) & X(I (Y Z)) 
\\
& Y(X(Z I)) \ar[dr]_{\beta_{Y,X,Z I}} \ar[r]^{11\rho_Z} & Y(XZ) \ar[ru]^{\beta_{Y,X,Z}} & X(Y(I Z)) \ar[ur]_{1\beta_{Y,I,Z}} 
\\
&& X(Y(Z I)) \ar[ru]_{1 1 c_{Z,I}} \ar[uur]^(.25){11\rho_Z} }
$$
Here we use the definition of $c$ and the coherence \eqref{ubc4}.
\epf

A b-functor $F:\C\to\D$ between unital b-categories is {\em unital} if it comes equipped with an isomorphism $\phi:F(I)\to I$ such that for every $X\in\C$ the following diagram
\beq\lb{ubf}\xymatrix{F(X\ot I) \ar[rr]^{F_{X,I}} \ar[d]_{F(\rho_X)} && F(X)\ot F(I) \ar[d]^{1\ot\phi} \\
F(X) && F(X)\ot I \ar[ll]_{\rho_{F(X)}} }\eeq
commutes. 
\ble\lb{ubf1}
Let $F:\C\to\D$ be a unital b-functor between unital b-categories. Then the following diagram commutes for all $X,Y\in\C$:
\beq\lb{pru}\xymatrix{F(X\ot Y) \ar[rr]^{F_{X,Y}} \ar[d]_{F(c_{X,Y})} && F(X)\ot F(Y) \ar[d]^{c_{F(X),F(Y)}} \\
F(Y\ot X) \ar[rr]^{F_{Y,X}} && F(Y)\ot F(X) }\eeq
\ele
\bpf
This follows from commutativity of the diagram
$$\xymatrix{ 
F(X (Y I))\ar[rr]^{F_{X,Y I}} \ar[ddd]_{F(\beta_{X,Y,I})} \ar[dr]^{F(1\rho_Y)} && F(X) F(Y I) \ar[rr]^{1F_{Y,I}} \ar[d]^{1F(\rho_Y)} && F(X) (F(Y) F(I)) \ar[ddd]^{\beta_{F(X),F(Y),F(I)}} \ar[ld]_{11\phi} \\ 
& F(XY) \ar[r]^{F_{X,Y}} \ar[d]_{F(c_{X,Y})} & F(X)F(Y) \ar[d]^{c_{F(X),F(Y)}} & F(X)(F(Y)I) \ar[l]_{1\rho_{F(Y)}} \ar[d]^{\beta_{F(X),F(Y),I}} 
\\ 
& F(XY) \ar[r]^{F_{X,Y}}  & F(X)F(Y)  & F(X)(F(Y)I) \ar[l]^{1\rho_{F(X)}}
\\
F(Y(X I)) \ar[rr]^{F_{Y,X I}} \ar[ur]_{F(1\rho_X)} && F(Y) F(X I) \ar[rr]^{1 F_{X,I}} \ar[u]_{1F(\rho_X)} && F(Y)(F(X) F(I)) \ar[ul]^{11\phi}
}$$
\epf

Here is the first main result of the paper:

\bth\lb{main}
(i) 
The structure of unital b-category is equivalent to the structure of braided monoidal category.
\nl
\smallskip\noindent (ii) 
The structure of symmetric unital b-category is equivalent to the structure of symmetric monoidal category.
\nl
\smallskip\noindent (iii) 
The structure of a unital b-functor between unital b-categories is equivalent to the structure of braided monoidal functor.
\eth
We prove it in the next three sections.

\bre 
Using Theorem \ref{main}, it is easy to see, for example, that $\beta_{X,I,Y}^{-1} = \beta_{I,X,Y}$ (since $c_{X,I}^{-1} = c_{I,X}$, see \cite[Prop.\,2.1]{js}). This is more cumbersome to check when starting directly from diagrams \eqref{cbc}, \eqref{ubc3} and \eqref{ubc4}.
\ere

\bre
The constructions of Sections \ref{bmbc} and \ref{bbfunct} below give a functor 
$${\bf BMCat}\ \to \ {\bf BCat}$$
from the 2-category of braided monoidal categories and braided monoidal functors to the 2-category of b-categories and b-functors.
The constructions of Sections \ref{ubcbmc} and \ref{bbfunct} below promote this functor to an equivalence
$${\bf BMCat}\ \to \ {\bf UBCat}$$
with target the 2-category of unital b-categories and unital b-functors.
\ere

\subsection{From braided monoidal categories to unital b-categories}\lb{bmbc}

Let $\C$ be a braided monoidal category with associativity $\alpha_{X,Y,Z}:X\ot(Y\ot Z)\to (X\ot Y)\ot Z$, braiding $c_{X,Y}:X\ot Y\to Y\ot X$, 
	and unit isomorphisms $\rho_X : X \otimes I \to X$ and $\lambda_X : I \otimes X \to X$.
Define the collection of isomorphisms
$\beta_{X,Y,Z}:X\ot(Y\ot Z)\to Y\ot(X\ot Z)$ as the composition
$$
\beta_{X,Y,Z} = \Big(
\xymatrix{X\ot(Y\ot Z)\ar[rr]^{\alpha_{X,Y,Z}} && (X\ot Y)\ot Z \ar[rr]^{c_{X,Y}\ot 1} && (Y\ot X)\ot Z \ar[rr]^{\alpha^{-1}_{Y,X,Z}}& & Y\ot(X\ot Z)} \Big) \ .
$$
Due to one of the hexagon axioms, $\beta_{X,Y,Z}$ coincides with the composition
$$
\beta_{X,Y,Z} = \Big(
\xymatrix{X\ot(Y\ot Z)\ar[rr]^{c_{X,YZ}} && (Y\ot Z)\ot X \ar[rr]^{\alpha^{-1}_{Y,Z,X}} && Y\ot (Z\ot X) \ar[rr]^{1\ot c_{Z,X}}& & Y\ot(X\ot Z)}
\Big) \ .$$

The coherence \eqref{cbc} for $\beta$ follows from commutativity of the diagram
$$\xymatrix{
	Y(X(ZW))
\ar[r]_{1\alpha} \ar@/^20pt/[rrrr]^{1\beta} \ar[dd]^\alpha & 
  Y((XZ)W)
\ar[rr]_{1(c1)} \ar[d]_\alpha && Y((ZX)W) \ar[d]^\alpha & Y(Z(XW))  \ar[l]^{1\alpha} \ar[dd]_\alpha \ar@/^40pt/[dddd]^{\beta} \\
& (Y(XZ))W \ar[rr]_{(1c)1} \ar[d]^{\alpha1} && (Y(ZX))W \ar[d]_{\alpha1} & \\
(YX)(ZW) \ar[r]^\alpha \ar[d]^{c1} & ((YX)Z)W \ar[d]^{(c1)1} && ((YZ)X)W \ar[d]^{(c1)1} & (YZ)(XW) \ar[l]_\alpha \ar[d]_{c1} \\
(XY)(ZW) \ar[r]_\alpha & ((XY)Z)W && ((ZY)X)W & (ZY)(XW) \ar[l]^\alpha\\
X(Y(ZW)) \ar[u]_\alpha \ar@/^40pt/[uuuu]^{\beta} \ar@/_40pt/[dddd]_{1\beta} \ar[d]^{1\alpha} & && & Z(Y(XW)) \ar[u]^\alpha \ar[d]_{1\alpha} \\
X((YZ)W) \ar[r]^\alpha \ar[d]^{1(c1)} & (X(YZ))W \ar[uu]_{\alpha1} \ar[d]_{(1c)1} \ar[rruuu]^{c1} && (Z(YX))W \ar[uu]^{\alpha1} & Z((YX)W) \ar[l]_\alpha \\
X((ZY)W) \ar[r]_\alpha & (X(ZY))W \ar[rruuu]_{c1} \ar[d]^{\alpha1} && (Z(XY))W \ar[u]^{(1c)1} \ar[d]_{\alpha1} & Z((XY)W) \ar[u]^{1(c1)} \ar[l]^\alpha \\
& ((XZ)Y)W \ar[rr]^{(c1)1} && ((ZX)Y)W \\
X(Z(YW)) \ar[r]^\alpha \ar[uu]_{1\alpha} \ar@/_20pt/[rrrr]_{\beta} & (XZ)(YW) \ar[rr]^{c1} \ar[u]^\alpha && (ZX)(YW) \ar[u]_\alpha & Z(X(YW)) \ar[l]_\alpha \ar[uu]^{1\alpha} \ar@/_40pt/[uuuu]_{1\beta}}$$
The b-structure we get is unital:
	we take the data $I,\rho$ for a unital b-category as in the braided category $\C$.
	The original braiding of $\C$ coincides with the natural collection defined by \eqref{br}:
$$\xymatrix{
(X\ot Y)\ot I \ar[d]_{c_{X,Y}\ot 1} \ar@/^20pt/[rrrr]^{\rho_{X\ot Y}} && X\ot(Y\ot I) \ar[ll]^{\alpha_{X,Y,I}} \ar[d]^{\beta_{X,Y,I}} \ar[rr]_{1\ot\rho_Y} && X\ot Y \ar[d]^{c_{X,Y}} \\ 
(Y\ot X)\ot I  \ar@/_20pt/[rrrr]_{\rho_{Y\ot X}} && Y\ot(X\ot I) \ar[ll]_{\alpha_{Y,X,I}} \ar[rr]^{1\ot\rho_X} && Y\ot X
}$$
The coherence \eqref{ubc4} for the unit isomorphism $\rho$ follows from commutative diagram
$$
\xygraph{!{0;/r6pc/:;/u4pc/::}
[]*+{Y(ZI)}
(
:[u]*+{Y(IZ)} ^{1c_{Z,I}}
 (
 :@/^20pt/[rrr]*+{I(YZ)}="u" ^{\beta_{Y,I,Z}}
 :[l]*+{(IY)Z}="e" ^{\alpha_{I,Y,Z}}
 ,
 :[r]*+{(YI)Z} _{\alpha_{Y,I,Z}}
 :"e" _{c_{Y,I}1}
 )
,
:[rrr]*+{(YZ)I} ^{\alpha_{Y,Z,I}}
 (
 :"u" _{c_{YZ,I}}
 ,
 :[l(1.5)d(.5)]*+{YZ}="b" ^{\rho_{YZ}}
 )
,
:"b" _{1\rho_Z}
)
}$$
Finally, commutativity of the coherence \eqref{ubc3} for the unital b-structure on $\C$ follows from the diagram:
$$\xygraph{!{0;/r8pc/:;/u8pc/::}
[]*+{(XY)(ZW)}
(
:[rr]*+{(ZW)(XY)}="u" ^{c_{XY,ZW}}
 (
 :[rd]*+{X((ZW)Y)}="ru" ^{\beta_{ZW,X,Y}}
  (
  :[d]*+{X(Y(ZW))}="s1" ^{1c_{ZW,Y}}
  ,
  :[l(.6)d(.3)]*+{(X(ZW))Y}="mr" _{\alpha_{X,ZW,Y}}
  :[l(.35)d(.675)]*+{((XZ)W)Y}="s2" _{\alpha_{X,Z,W}1}
  )
 ,
 :[l(.3)d(.6)]*+{((ZW)X)Y}="mu" ^{\alpha_{ZW,X,Y}}
 :"mr" ^{c_{ZW,X}1}
 )
,
:[ld]*+{Z((XY)W)} _{\beta_{X Y,Z,W}}
:[d]*+{Z(W(XY))} _{1c_{XY,W}}
 (
 :@/^30pt/"u" ^{\alpha_{Z,W,XY}}
 ,
 :[r(.675)u(.35)]*+{Z((WX)Y)} _{1\alpha_{W,X,Y}}
  (
  :[r(1.35)u(.7)]*+{(Z(WX))Y} ^{\alpha_{Z,WX,Y}}
   (
   :"mu" ^{\alpha_{Z,W,X}1}
   ,
   :[l(.3)d(.6)]*+{(Z(XW))Y}="m" ^{1c_{W,X}1}
   :[r(.75)d(.75)]*+{((ZX)W)Y}="mb" ^{\alpha_{Z,X,W}1}
   :"s2" ^{c_{Z,X}11}
   )
  ,
  :[r(.5)d(.5)]*+{Z((XW)Y)}="ml" ^{1c_{W,X}1}
  :"m" _{\alpha_{Z,XW,Y}}
  )
 ,
 :[rd]*+{Z(X(WY))} _{1\beta_{W,X,Y}}
  (
  :"ml" ^{1\alpha_{X,W,Y}}
  ,
  :[r(.8)u(.2)]*+{(ZX)(WY)} ^{\alpha_{Z,X,WY}}
   (
   :"mb" ^{\alpha_{ZX,W,Y}}
   ,
   :[r(.5)d(.5)]*+{(XZ)(WY)}="b" ^{c_{X,Z}1}
   :"s2" ^{\alpha_{XZ,W,Y}}
   )
  ,
  :[rd]*+{X(Z(WY))} _{\beta_{Z,X,WY}}
   (
   :"b" ^{\alpha_{X,Z,WY}}
   ,
   :@/_10pt/"ru" _(.6){1\alpha_{Z,W,Y}}
   ,
   :[ru]*+{X(Z(YW))} _{11c_{W,Y}}
   :"s1" _{1\beta_{Z,Y,W}}
   )
  )
 )
)
}$$

If $\C$ is a symmetric monoidal category, the diagram
$$
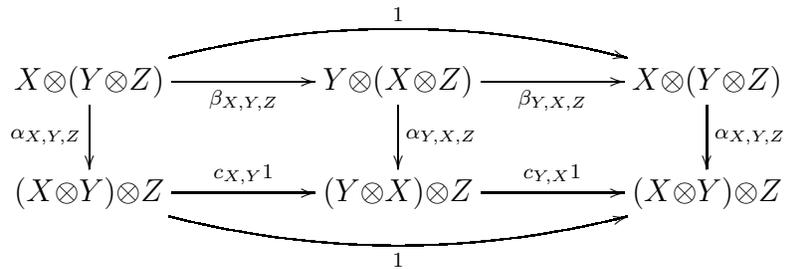
\xymatrix{
X\ot (Y\ot Z) \ar[rr]_{\beta_{X,Y,Z}} \ar@/^20pt/[rrrr]^1 \ar[d]_{\alpha_{X,Y,Z}} && Y\ot (X\ot Z) \ar[rr]_{\beta_{Y,X,Z}} \ar[d]^{\alpha_{Y,X,Z}} &&  X\ot (Y\ot Z) \ar[d]^{\alpha_{X,Y,Z}}\\
(X\ot Y)\ot Z \ar[rr]^{c_{X,Y}1} \ar@/_20pt/[rrrr]_1  && (Y\ot X)\ot Z \ar[rr]^{c_{Y,X}1}  &&  (X\ot Y)\ot Z 
}$$
shows that it is symmetric as a b-category.

\subsection{From unital b-categories to braided monoidal categories}\lb{ubcbmc}

Let now $\C$ be a unital b-category. Define the collection of isomorphisms $c_{X,Y}:X\ot Y\to Y\ot X$ as in \eqref{br} and $\alpha_{X,Y,Z}:X\ot(Y\ot Z)\to (X\ot Y)\ot Z$ as 
$$
\alpha_{X,Y,Z} = \Big(
\xymatrix{X\ot (Y\ot Z) \ar[rr]^{1\ot c^{-1}_{Z,Y}} && X\ot (Z\ot Y) \ar[rr]^{\beta^{-1}_{Z,X,Y}} && Z\ot (X\ot Y) \ar[rr]^{c_{Z,X\ot Y}}& & (X\ot Y)\ot Z}
\Big) \ .
$$
Note that the diagram \eqref{ubc2} implies that $\alpha_{X,Y,Z}$ coincides with the composition
$$
\alpha_{X,Y,Z} = \Big(
\xymatrix{X\ot (Y\ot Z) \ar[rr]^{1\ot c_{Y,Z}} && X\ot (Z\ot Y) \ar[rr]^{\beta_{X,Z,Y}} && Z\ot (X\ot Y) \ar[rr]^{c^{-1}_{X\ot Y,Z}}& & (X\ot Y)\ot Z} \Big) \ .
$$
The hexagon axioms follow from the commutative diagrams
$$\xymatrix{
& (X\ot Y)\ot Z \ar[rr]^{c_{X\ot Y,Z}} && Z\ot(X\ot Y) \ar[rdd]^{\alpha_{Z,X,Y}} \\
& X\ot(Z\ot Y) \ar[rru]^{\beta_{X,Z,Y}} && Z\ot(Y\ot X) \ar[u]^{1\ot c_{Y,X}} \\
X\ot(Y\ot Z) \ar[ruu]^{\alpha_{X,Y,Z}} \ar[ru]_{1\ot c_{Y,Z}} \ar[rdd]_{1\ot c_{Y,Z}} &&&& (Z\ot X)\ot Y \\
& Y\ot(X\ot Z) \ar[lu]_{\beta_{Y,X,Z}} \ar[rr]^{1\ot c_{X,Z}} \ar[rrd]^{c_{Y,X\ot Z}} && Y\ot(Z\ot X) \ar[uu]^{\beta_{Y,Z,X}} \ar[ru]^{c_{Y,Z\ot X}} \\
& X\ot(Z\ot Y) \ar[rr]_{\alpha_{X,Z,Y}} && (X\ot Z)\ot Y \ar[ruu]_{c_{X,Z}\ot 1} 
}$$
and
$$\xymatrix{
& X\ot (Y\ot Z) \ar[rr]^{c_{X,Y\ot Z}} \ar[ldd]_{\alpha_{X,Y,Z}} \ar[d]^{1\ot c_{Y,Z}} \ar[drr]^{\beta_{X,Y,Z}} && (Y\ot Z)\ot X  \\
& X\ot(Z\ot Y) \ar[dd]^{\beta_{X,Z,Y}} && Y\ot(X\ot Z) \ar[dr]_{1\ot c_{X,Z}} \\
(X\ot Y)\ot Z \ar[rdd]_{c_{X,Y}\ot 1} \ar[rd]^{c_{X\ot Y,Z}}  &&&& Y\ot(Z\ot X) \ar[ld]_{\beta_{Y,Z,X}} \ar[luu]_{\alpha_{Y,Z,X}} \\
& Z\ot(X\ot Y)  \ar[rr]^{1\ot c_{X,Y}}  && Z\ot(Y\ot X)  \\
& (Y\ot X)\ot Z \ar[rru]_{c_{Y\ot X,Z}} && Y\ot(X\ot Z) \ar[ruu]_{1\ot c_{X,Z}} \ar[ll]^{\alpha_{Y,X,Z}}
}$$
	In both cases the middle hexagon commutes by Lemma \ref{brc}.

The pentagon axiom follows from the commutative diagram
$$\xygraph{!{0;/r3.5pc/:;/u3.5pc/::}
[]*+{(ZW)(XY)}
(
:[u(1.1)r(1.7)]*+{(XY)(ZW)}="t" ^{c_{ZW,XY}}
 (
 :[d(1.1)r(1.7)]*+{(XY)(WZ)} _{1c_{Z,W}}
 :[d(1.1)r(1.7)]*+{W((XY)Z)}="tr" _{\beta_{XY,W,Z}}
 :[l(.6)d(2)]*+{W(Z(XY))}="e" _{1c_{XY,Z}}
 ,
 :[d(1.65)r(.1)]*+{Z((XY)W)} _{\beta_{XY,Z,W}}
 :[d(1.65)r(.1)]*+{Z(W(XY))}="m" _{1c_{XY,W}}
 :"e" ^{\beta_{Z,W,XY}}
 ,
 :@/^30pt/[d(3.3)r(5.1)]*+{((XY)Z)W}="r" ^{\alpha_{XY,Z,W}}
  (
  :"tr" ^{c_{(XY)Z,W}}
  ,
  :[l(.6)d(2)]*+{(Z(XY))W}="mr" _{c_{XY,Z}1}
  :"e" ^{c_{Z(XY),W}}
  )
 )
,
:[d(1.1)l(1.7)]*+{X((ZW)Y)} _{\beta_{ZW,X,Y}}
:[d(1.1)l(1.7)]*+{X(Y(ZW))} ^{1c_{ZW,Y}}
 (
 :@/^30pt/"t" ^{\alpha_{X,Y,ZW}}
 ,
 :@/_30pt/[d(6)r(1.8)]*+{X((YZ)W)} _{1\alpha_{Y,Z,W}}
  (
  :@/_30pt/[r(6.6)]*+{(X(YZ))W} _{\alpha_{X,YZ,W}}
   (
   :@/_30pt/"r" _{\alpha_{X,Y,Z}1}
   ,
   :[l(2.2)]*+{W(X(YZ))}="b" _{c_{X(YZ),W}}
   :[u(2)r(.6)]*+{W(X(ZY))}="br" ^{11c_{Y,Z}}
   :"e" _{1\beta_{X,Z,Y}}
   ,
   :[u(2)r(.6)]*+{(X(ZY))W} ^{1c_{Y,Z}1}
    (
    :"mr" ^{\beta_{X,Z,Y}1}
    ,
    :"br" _{c_{X(ZY),W}}
    )
   )
  ,
  :[u(2)l(.6)]*+{X(W(YZ))}="l" _{1c_{YZ,W}}
   (
   :"b" _{\beta_{X,W,YZ}}
   ,
   :[r(2.2)]*+{X(W(ZY))}="bl" _{11c_{Y,Z}}
   :"br" _{\beta_{X,W,ZY}}
   )
  )
 ,
 :[d(2)r(.6)]*+{X(Y(WZ))} ^{11c_{Z,W}}
 :"l" ^{1\beta_{Y,W,Z}}
 ,
 :[r(1.4)d(1.2)]*+{X(Z(YW))} ^{1\beta_{Y,Z,W}}
  (
  :[r(1.4)d(1.2)]*+{X(Z(WY))} ^{11c_{Y,W}}
   (
   :[r(1.4)u(1.2)]*+{Z(X(WY))}="ml" _{\beta_{X,Z,WY}}
   :"m" _{1\beta_{X,W,Y}}
   ,
   :"bl" ^{1\beta_{Z,W,Y}}
   )
  ,
  :[r(1.4)u(1.2)]*+{Z(X(YW))} _{\beta_{X,Z,YW}}
  :"ml" ^{11c_{Y,W}}
  )
 )
)
}$$

Define the collection of isomorphisms
$\lambda_X:I\ot X\to X$ natural in $X\in\C$ as the compositions
$$\xymatrix{I\ot X\ar[rr]^{c_{I,X}} && X\ot I\ar[rr]^{\rho_X} && X}$$
The axiom for the monoidal unit object follows from the commutative diagram
$$\xymatrix{
X\ot(I\ot Y) \ar[rrrrrr]^{\alpha_{X,I,Y}} \ar[drr]^{1\ot c_{I,Y}} \ar@/_10pt/[ddrr]_{1\ot\lambda_Y} &&&&&& (X\ot I)\ot Y \ar[dll]_{c_{X\ot I,Y}} \ar@/^50pt/[ddllll]^(.34){\rho_X\ot 1} \\ 
&& X\ot(Y\ot I) \ar[rr]^{\beta_{X,Y,I}} \ar[d]^{1\ot\rho_Y} && Y\ot(X\ot I) \ar[d]^{1\ot\rho_X} \\ 
&& X\ot Y \ar[rr]^{c_{X,Y}} && Y\ot X}$$

Finally if $\C$ is a symmetric b-category the diagram
$$\xymatrix{X\ot Y \ar[rr]_{c_{X,Y}} \ar@/^20pt/[rrrr]^1 && Y\ot X \ar[rr]_{c_{Y,X}} &&  X\ot Y \\
X\ot (Y\ot I) \ar[rr]^{\beta_{X,Y,I}} \ar@/_20pt/[rrrr]_1 \ar[u]^{1\rho_Y} && Y\ot (X\ot I) \ar[rr]^{\beta_{Y,X,I}} \ar[u]_{1\rho_X} &&  X\ot (Y\ot I) \ar[u]_{1\rho_Y}
}$$
shows that it is symmetric as a monoidal category.

\subsection{Unital b-functors vs braided monoidal functors}\lb{bbfunct}

Here we prove the last part of the Theorem \ref{main}. Namely for a braided monoidal functor we define its unital b-structure and visa versa.
In both cases the structural isomorphisms $F_{X,Y} : F(X\ot Y)\to F(X)\ot F(Y)$ and $\phi:F(I)\to I$ are the same. 
We simply prove that one set of coherence conditions is equivalent to the other. 

\medskip

Let $F:\C\to\D$ be a braided monoidal functor with monoidal isomorphism $F_{X,Y}:F(X\ot Y)\to F(X)\ot F(Y)$.
The following diagram shows that this isomorphism is a b-functor structure, i.e.\ the diagram \eqref{cbf} commutes for $F_{X,Y}$.
$$\xymatrix{ 
F(X (Y Z))\ar[rr]^{F_{X,Y Z}} \ar[ddd]_{F(\beta_{X,Y,Z})} \ar[dr]^{F(\alpha_{X,Y,Z})} && F(X) F(Y Z) \ar[rr]^{1F_{Y,Z}}  && F(X) (F(Y) F(Z)) \ar[ddd]^{\beta_{F(X),F(Y),F(Z)}} \ar[ld]_{\alpha_{F(X),F(Y),F(Z)}} \\ 
& F((XY)Z) \ar[r]^{F_{XY,Z}} \ar[d]_{F(c_{X,Y}1)} & F(XY)F(Z) \ar[r]^(.4){F_{X,Y}1} \ar[d]^{F(c_{X,Y})1} & (F(X)F(Y))F(Z) \ar[d]^{c_{F(X),F(Y)}1} 
\\ 
& F((YX)Z) \ar[r]^{F_{YX,Z}}  & F(YX)F(Z) \ar[r]^(.4){F_{Y,X}1} & (F(Y)F(X))F(Z) 
\\
F(Y(X Z)) \ar[rr]_{F_{Y,X Z}} \ar[ur]_{F(\alpha_{Y,X,Z})} && F(Y) F(X Z) \ar[rr]_{1 F_{X,Z}}  && F(Y)(F(X) F(Z)) \ar[ul]^{\alpha_{F(Y),F(X),F(Z)}}
}$$
Unitality of $F$ is automatic.

Now let $F:\C\to\D$ be a unital b-functor.
The braiding axiom for $F$ follows from Lemma \ref{ubf1}.
The monoidal coherence for $F$ comes from commutativity of the diagram:
$$\xymatrix{ 
F(X (Y Z))\ar[rr]^{F_{X,Y Z}} \ar[ddd]_{F(\alpha_{X,Y,Z})} \ar[dr]^{F(1c_{Y,Z})} && F(X) F(Y Z) \ar[rr]^{1F_{Y,Z}} \ar[d]^{1F(c_{Y,Z})} && F(X) (F(Y) F(Z)) \ar[ddd]^{\alpha_{F(X),F(Y),F(Z)}} \ar[ld]_{1c_{F(Y),F(Z)}} \\ 
& F(X(ZY)) \ar[r]^{F_{X,ZY}} \ar[d]_{F(\beta_{X,Z,Y})} & F(X)F(ZY) \ar[r]^{1F_{Z,Y}} & F(X)(F(Z)F(Y)) \ar[d]^{\beta_{F(X),F(Z),F(Y)}} 
\\ 
& F(Z(XY)) \ar[r]^{F_{Z,XY}}  & F(Z)F(XY) \ar[r]^{1F_{X,Y}} & F(Z)(F(X)F(Y)) 
\\
F((XY)Z) \ar[rr]^{F_{XY, Z}} \ar[ur]_{c_{XY,Z}} && F(XY)F(Z) \ar[rr]^{F_{X,Y}1} \ar[u]_{c_{F(XY),F(Z)}} && (F(X)F(Y))F(Z) \ar[ul]^{c_{F(X)F(Y),F(Z)}}
}$$
The coherence \eqref{ubf} gives one of the unit preservation axioms (the right one). The left unit preservation axiom follows from
$$\xymatrix{
F(I\ot X) \ar[rrr]^{F_{I,X}} \ar[rd]^{F(c_{I,X})} \ar[ddd]_{F(\lambda_X)} &&& F(I)\ot F(X) \ar[dl]_{c_{F(I),F(X)} }\ar[ddd]^{\phi 1} \\
& F(X\ot I) \ar[r]^{F_{X,I}} \ar[ldd]_{F(\rho_X)} & F(X)\ot F(I) \ar[d]^{1\phi} \\
&& F(X)\ot I \ar[lld]^{\rho_{F(X)}} \\
F(X) &&& I\ot F(X) \ar[lll]^{\lambda_{F(X)}} \ar[ul]_{c_{I,F(X)}}
}$$

Thus the proof of Theorem \ref{main} is complete. 

\subsection{Pre-unital b-categories}

Here we define unital b-categories without the unit. This modification will be useful for treating some of the examples.

We call a b-category $\C$ {\em pre-unital} if it comes equipped with a natural (in $X,Y\in\C$) collection  $c_{X,Y}:X\ot Y\to Y\ot X$ of isomorphisms satisfying the coherence axioms \eqref{ubc3}, \eqref{brc1}, \eqref{ubc2}.

\bpr
A pre-unital b-category can be fully embedded into a unital b-category. 
\epr
\bpf
Let $\I$ be the one object one morphism category.
Consider the disjoint union $\tilde\C = \C\cup\I$.
Extend the tensor product from $\C$ to $\tilde\C$ by $X\ot I = X = I\ot X$.
Extend the b-structure from $\C$ to $\tilde\C$ by 
$$\tilde\beta_{I,Y,Z} = 1_{Y\ot Z},\qquad \tilde\beta_{X,I,Z} = 1_{X\ot Z},\qquad \tilde\beta_{X,Y,I} = c_{X,Y}\ .$$
To see that $\tilde\C$ is a b-category we need to verify the b-axiom for the extended b-structure, i.e.\ when some of the objects $X,Y,Z,W$ in the diagram \eqref{cbc} are $I$. Of all the cases the only non-trivial is when $W=I$.
In this case the diagram \eqref{cbc} is just \eqref{brc1}.

Finally define $\rho_X$ to be the identity $1_X$.
It is then straightforward (and tedious) to check that $\tilde\C$ is a unital b-category.
\epf

A b-functor $F:\C\to \D$ between pre-unital b-categories is {\em pre-unital} if it satisfies the coherence \eqref{pru}.
It is straightforward to see that a pre-unital b-functor $F:\C\to \D$ between pre-unital b-categories extends to a 
unital b-functor $\tilde F:\tilde\C\to \tilde\D$ between unital b-categories via 
	$\tilde F(I) = I$,  $\tilde F_{I,X} = \tilde F_{X,I} = 1_{F(X)}$ and $\phi = 1_I$.

\section{Examples of b-categories}\label{sec:examples}

Let $k$ be a field. 
Here we look at
  (pre-unital) b-structures on the category $\Vect_k$ 
of (finite dimensional) vector spaces over $k$. Denote by $\ot=\ot_k$ the standard tensor product of vector spaces.
Note that any $k$-linear tensor product on $\Vect_k$ necessarily has the form
$$\ot_M:\Vect_k\times\Vect_k\to\Vect_k,\qquad (U,V)\mapsto U\ot_MV = U\ot M\ot V$$
for some $M\in\Vect_k$. Indeed, a $k$-linear tensor product on $\Vect_k$ is determined by its value (in this case $M$) on the pair $(k,k)$ of one-dimensional vector spaces.

\medskip

Below, we first recall from \cite{da} the description of semi-groupal (monoidal without the unit object) structures on $\Vect_k,\ot_M$, and then turn to the description of b-structures. 

\medskip\noindent
{\bf Notation:}
For an operator $P$ on $M^{\ot 2}$, the operator $P_{ij}$ on $M^{\ot n}$ is $P$ acting on the $i$-th and $j$-th components. 

\subsection{Semi-groupal structures on the category of vector spaces}

A category $\C$ with a tensor product $\ot:\C\times\C\to\C$ is {\em semi-groupal} if it is equipped with a natural in $X,Y,Z\in\C$ collection of isomorphisms $\alpha_{X,Y,Z}:X\ot(Y\ot Z)\to (X\ot Y)\ot Z$ (an {\em associativity constraint}) satisfying the pentagon axiom.

\begin{prop}[{\cite[Prop.\,7.1]{da}}]
An associativity constraint for $\Vect_k,\ot_M$ corresponds to a solution $\Phi\in Aut_k(M^{\ot 2})$ of the {\em pentagon equation}:
$$\Phi_{12}\Phi_{13}\Phi_{23} = \Phi_{23}\Phi_{12}.$$
\epr

\bre
Explicitly the associativity for $\Vect_k,\ot_M$ gives an automorphism $\Phi\in Aut_k(M^{\otimes
2})$: 
$$\xymatrix{M\ot M = k\ot_M (k\ot_M k)\ar[rr]^{\alpha_{k,k,k}} && (k\ot_M k)\ot_M k = M\ot M}$$ and has the
following general form:
$$\xymatrix{ 
U_1\ot_M(U_2\ot_M U_3) \ar@{=}[d] \ar[rr]^{\alpha_{U_1, U_2 ,U_3}} && (U_1\ot_M U_2)\ot_M U_3 \ar@{=}[d] \\ 
U_1\ot M\ot U_2\ot M\ot U_3  \ar[rr]^{\Phi_{24}} && U_1\ot M\ot U_2\ot M\ot U_3}$$ 

\ere

\subsection{b-structures on vector spaces and the Yang-Baxter equation}

\bpr \label{b-cat-on-vect}
A b-structure on $\Vect_k,\ot_M$ corresponds to a solution $B\in Aut_k(M^{\ot 2})$ of the {\em hexagon equation}:
$$B_{12}B_{23}B_{12} = B_{23}B_{12}B_{23}\ .$$
\epr
\bpf
By naturality a b-structure $\beta$ for $\Vect_k,\ot_M$ is given by an automorphism $B\in Aut_k(M^{\otimes
2})$: 
$$\xymatrix{M\ot M = k\ot_M (k\ot_M k)\ar[rr]^{\beta_{k,k,k}} && (k\ot_M k)\ot_M k = M\ot M}$$ and has the
following general form 
$$\xymatrix{ 
U_1\ot_M(U_2\ot_M U_3) \ar@{=}[d] \ar[rr]^{\beta_{U_1, U_2 ,U_3}} && U_2\ot_M (U_1\ot_M U_3) \ar@{=}[d] \\ 
U_1\ot M\ot U_2\ot M\ot U_3  \ar[rr]^{B_{24}t_{13}} && U_2\ot M\ot U_1\ot M\ot U_3}$$ 
Here $t:U_1\ot U_2\to U_2\ot U_1$ is the ordinary interchange of tensor factors.
In particular, 
$$\beta_{k,k,k\ot_M k} = \beta_{k,k,M} = B_{12}\in Aut_k(M^{\otimes 3}).$$ 
Combining this with 
$$1_k\ot_M\beta_{k,k,k} = B_{23}\in Aut_k(M^{\otimes 3}) \ ,$$ 
we see that the coherence \eqref{cbc} for $\beta$ is equivalent to
the hexagon equation for $B$. Indeed, by naturality the general coherence is equivalent to the specialisation of the coherence diagram \eqref{cbc} to $X,Y,Z,W=k$, which in the equational form is
$$\beta_{k,k,M}\circ(1\ot\beta_{k,k,k})\circ \beta_{k,k,M} = (1\ot\beta_{k,k,k})\circ\beta_{k,k,M}\circ(1\ot\beta_{k,k,k}) \ .$$
\epf

A $k$-linear autoequivalence of $\Vect_k$ has to be the identity.
A b-functor structure on $Id:\Vect_k\to\Vect_k$ amounts to an automorphism $g:M\to M$ such that
$$(g\ot g)B = B(g\ot g)\ .$$
Composition of b-functors corresponds to composition of automorphisms $M\to M$.

\subsection{Pre-unital b-structures on the category of vector spaces}

\bpr 
A structure of a pre-unital b-category on $\Vect_k,\ot_M$ corresponds to a solution $B\in Aut_k(M^{\ot 2})$ of the hexagon equation together with 
   a solution $C\in Aut_k(M)$ of
\beq\lb{mfu}B(1\ot C)B = (1\ot C)B(1\ot C),\qquad B(1\ot C^2)B = (C^2\ot 1)\eeq and
\beq\lb{mfu1}B_{23}C_3B_{12}B_{23}C_2B_{32} = C_1C_2B_{32}t_{23} \ . \eeq
\epr
\bpf
Let $(\beta, c)$ be a structure of pre-unital b-category on $\Vect_k,\ot_M$.
By naturality, the collection $c$ is given by an automorphism $C\in Aut_k(M)$: 
$$
	C = \Big(
\xymatrix{M = k\ot_M k\ar[rr]^{c_{k,k}} && k\ot_M k = M} \Big) \ ,$$ and has the
following general form 
$$\xymatrix{ 
U_1\ot_M U_2 \ar@{=}[d] \ar[rr]^{c_{U_1, U_2}} && U_2\ot_M U_1 \ar@{=}[d] \\ 
U_1\ot M\ot U_2  \ar[rr]^{(1\ot C\ot 1)t_{13}} && U_2\ot M\ot U_1}$$
In particular, 
$$c_{k,k\ot_M k} = c_{k,M} = (1\ot C)t,\quad c_{k\ot_M k,k} = c_{M,k} = (C\ot 1)t\quad \in Aut_k(M^{\otimes 2}).$$ 
Combining it with 
$$1_k\ot_Mc_{k,k} = (1\ot C)\in Aut_k(M^{\otimes 2})$$ 
	we see that the coherences \eqref{brc1} and \eqref{ubc2} for $\beta$ and $c$ are equivalent to the equations \eqref{mfu}. 
Indeed, by naturality the general coherences are equivalent to the specialisations of the coherence diagrams \eqref{brc1}, \eqref{ubc2} to $X,Y,Z,W=k$, which in the equational form are
$$\beta_{k,k,k}\circ(1\ot c_{k,k})\circ \beta_{k,k,k} = (1\ot c_{k,k})\circ\beta_{k,k,k}\circ(1\ot c_{k,k}),$$
$$\beta_{k,k,k}\circ(1\ot c_{k,k})^2\circ \beta_{k,k,k} = c_{M,k}\circ c_{k,M}.$$
Finally the coherence \eqref{ubc3} is equivalent to the equation
$$(1\ot c_{M,k})\circ\beta_{M,k,k}\circ c_{M,M} = (1\ot\beta_{k,k,k})\circ(1\ot1\ot c_{k,k})\circ\beta_{k,k,M}\circ(1\ot\beta_{k,k,k})\circ(1\ot c_{M,k})\circ\beta_{M,k,k}.$$
Using that $c_{M,M} = C_2t_{13}, \beta_{M,k,k} = B_{23}t_{12}\in Aut_k(M^{\otimes 3})$ we get the condition \eqref{mfu1}. 
\epf

\subsection{Categorical b-magmas}\lb{pbc}

Let us call a b-category $\C$ {\em categorical b-magma} if $\C$ is a 
groupoid and the tensor product maps on automorphisms
$$\C(X,X)\to \C(X\ot Y,X\ot Y)\leftarrow \C(Y,Y)$$
are isomorphisms for all 
$X,Y\in\C$. The last condition implies that the automorphism groups $\C(X,X)$ are isomorphic for all $X\in\C$ and must be abelian (by the Eckmann-Hilton argument). 
Clearly, the set $\pi_0(\C)$ of isomorphisms classes of objects of a categorical b-magma $\C$ is a b-magma with the binary operation $[X][Y] = [X\ot Y]$. Denote by $\pi_1(\C)$ the automorphism group $\C(X,X)$ of an object $X\in\C$.

Choose representatives $s:\pi_0(\C)\to \C$ of the isomorphism classes of objects, as well as isomorphisms $t_{x,y} : s(x) \ot s(y) \to s(xy)$ for each pair $x,y \in \pi_0(\C)$. 
The coherence isomorphisms of a b-category evaluated on objects give rise to non-zero scalars $r(x,y,z)\in \pi_1(\C)$ via
\begin{align*}
r(x,y,z) \, 1_{s(x(yz))}
=
\Big( &
s(x(yz)) 
\xrightarrow{t_{x,yz}^{-1}} s(x)s(yz) 
\xrightarrow{1\,t_{y,z}^{-1}} s(x)(s(y)s(z)) 
\\
&
\xrightarrow{\beta_{s(x),s(y),s(z)}} s(y)(s(x)s(z)) 
\xrightarrow{1\,t_{x,z}} s(y)(s(xz)) 
\xrightarrow{t_{y,xz}} s(y(xz)) = s(x(yz)) 
\Big) \ .
\end{align*}
The coherence axiom for the isomorphisms evaluated on objects amounts to the equation
$$r(y,z,xw)r(x,z,w)r(x,y,zw) = r(x,y,w)r(x,z,yw)r(y,z,w)$$ for all $x,y,z,w\in \pi_0(\C)$. 

A different choice of section $s' : \pi_0(\C)\to \C$ and isomorphisms $t'_{x,y} : s'(x) \ot s'(y) \to s'(xy)$ produces in general a different set of values $r'_{x,y,z}$. To relate $r'$ to $b$, pick isomorphisms $u_x : s(x) \to s'(x)$. These determine constants $q_{x,y} \in \pi_1(\C)$ ($x,y \in \pi_0(\C)$) via
$q_{x,y} \cdot u_{xy} \circ t_{x,y} = t'_{x,y} \circ (u_x \ot u_y)$.
One computes that, independent of the choice of $u$'s, for all $x,y,z\in \pi_0(\C)$,
$$
  r'(x,y,z) = q^{\phantom{-1}}_{x,z} \, q^{-1}_{x,yz} \, q^{\phantom{-1}}_{y,xz} \, q^{-1}_{y,z} \, r(x,y,z) \ .
$$

A b-functor $F:\C\to \C'$ between categorical b-magmas gives rise to a morphism $f:\pi_0(\C)\to \pi_0(\C')$ of b-magmas. Assume we have chosen sections $s,s'$ and isomorphisms $t,t'$ for $\C$ and $\C'$, and that in addition we have fixed isomorphisms $h_x : s'(f(x)) \to F(s(x))$.
The coherence isomorphisms of $F:\C\to \C'$ evaluated on objects determine a collection of constants $q(x,y) \in  \pi_1(\C')$ ($x,y \in \pi_0(\C)$) via
\begin{align*}
 q(x,y) \, 1_{s'(f(xy))} = \Big( & 
 s'(f(xy)) \xrightarrow{h_{xy}} F(s(xy))
 \xrightarrow{F(t_{x,y}^{-1})}  F(s(x) \otimes s(y))
 \xrightarrow{F_{s(x),s(y)}} F(s(x)) \otimes F(s(y)) \\
 & \xrightarrow{h_x^{-1} \ot h_y^{-1}} s'(f(x)) \otimes s'(f(y)) 
 \xrightarrow{{t'}_{f(x),f(y)}}  s'(f(x)f(y)) = s'(f(xy)) 
 \Big) \ .
\end{align*} 
The coherence axiom \eqref{cbf} for a b-functor evaluated on objects now amounts to the equation
$$r'(f(x),f(y),f(z)) \, q(y,z) \, q(x,yz) = q(x,z) \, q(y,xz) \, r(x,y,z)$$ for all $x,y,z\in \pi_0(\C)$. 

Let $\tilde F : \C \to \C'$ be another b-functor as above and let $a :  F \Rightarrow \tilde F$ be a b-transformation. Since $a$ is invertible $\tilde F$ induces the same map $f : \pi_0(\C) \to \pi_0(\C')$ as $F$. The b-transformation $a$ gives rise to constants $p(x)$ for all $x \in \pi_0(\C)$ via
$$
  p(x) = \Big( s'(f(x)) \xrightarrow{h_x} F(s(x)) \xrightarrow{a_{s(x)}} \tilde F(s(x)) \xrightarrow{\tilde h^{-1}_x} s'(f(x)) \Big) \ .
$$
The coherence condition \eqref{cbt} for a b-transformation becomes $q(x,y)p(x)p(y)=p(x,y)\tilde q(x,y)$.

\section{b-bicategories} \label{sec:b-bicategories}

Here we define a bicategorical analogue of b-category.
By a {\em tensor product} on a bicategory $\bC$ we will mean a functor 
$$\bC\times\bC\to\bC,\qquad (X,Y)\mapsto X\ot Y\ .$$
In particular for any 1-morphisms $f:X\to Y,\ g:Z\to W$ we a have an invertible 2-cell
\beq\lb{sli}\xygraph{ !{0;/r6.5pc/:;/u6.5pc/::} []*+{X\ot Z}
(
:[r]*+{X\ot W} ^{1\ot g}
:[d]*+{Y\ot W}="e" ^{f\ot 1}
,
:[d]*+{Y\ot Z} _{f\ot 1}
:"e" _{1\ot g}
,
:@{}[r(.6)d(.5)]*+{\Downarrow \scriptstyle{s_{f,g}}}
)
}\eeq
with bi-multiplicativity property for the collection $s$. 

\subsection{Definition}

A bicategory $\bC$ is a {\em b-bicategory} if it is equipped with a  tensor product together with a pseudo-natural in $X,Y,Z\in\bC$ collection of 1-equivalences
$$\beta_{X,Y,Z}:X\ot(Y\ot Z)\to Y\ot(X\ot Z)\ ,$$
and natural in $X,Y,Z,W\in\bC$ collection of 2-isomorphisms
\beq\lb{2cbc}\xymatrix{ & Y\ot(X\ot (Z\ot W)) \ar[rr]^{1\ot\beta_{X,Z,W}} && Y\ot(Z\ot (X\ot W))  \ar[rd]^{\ \beta_{Y,Z,X\ot W}}\\
X\ot(Y\ot (Z\ot W)) \ar[ru]^{\beta_{X,Y,Z\ot W}} \ar[rd]_{1\ot\beta_{Y,Z,W}} &&\ \Downarrow\scriptstyle{\gamma_{X,Y,Z,W}} && Z\ot(Y\ot (X\ot W))\\
& X\ot(Z\ot (Y\ot W)) \ar[rr]^{\beta_{X,Z,Y\ot W}} && Z\ot(X\ot (Y\ot W)) \ar[ru]_{\ 1\ot\beta_{X,Y,W}}} \ .\eeq
	Pseudo-naturality of $\beta$ means that for all 1-morphisms $f:X\to Y,\ g:Z\to W,\ h:U\to V$ we have an invertible 2-cell
\beq\lb{nat}\xygraph{ !{0;/r8.5pc/:;/u6.5pc/::} []*+{X\ot (Z\ot U)}
(
:[r]*+{Z\ot(X\ot U)} ^{\beta_{X,Z,U}}
:[d]*+{W\ot(Y\ot V)}="e" ^{f\ot g\ot h}
,
:[d]*+{Y\ot(W\ot V)} _{f\ot g\ot h}
:"e" _{\beta_{Y,W,V}}
,
:@{}[r(.6)d(.5)]*+{\Downarrow \scriptstyle{b_{f,g,h}}}
)
}\eeq
with tri-multiplicativity property for the collection $b$. 
The 2-isomorphisms \eqref{2cbc} and \eqref{nat} must be such that the following two pastings of 2-cells 
agree
	(the squares in the pasting schemes below are filled with the 2-cells $b_{1,1,\beta}$):

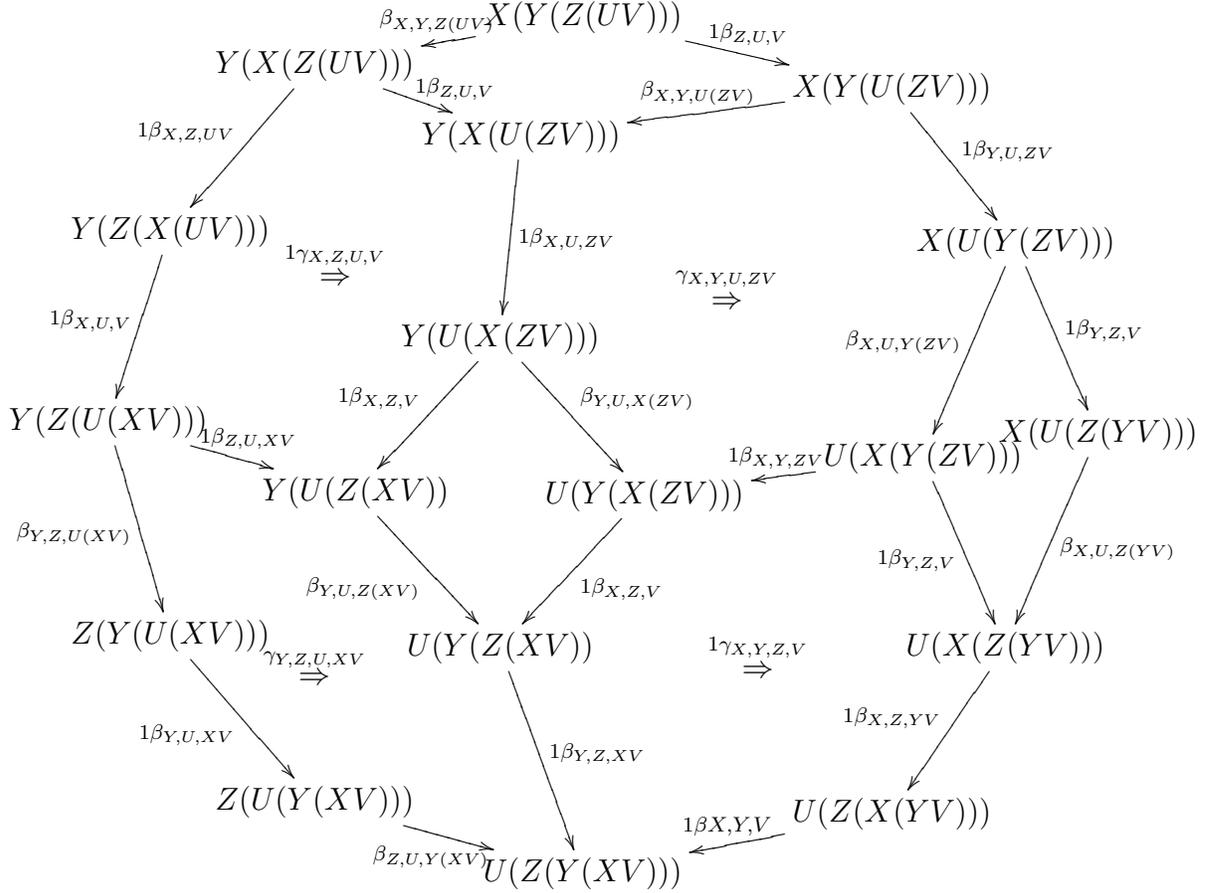
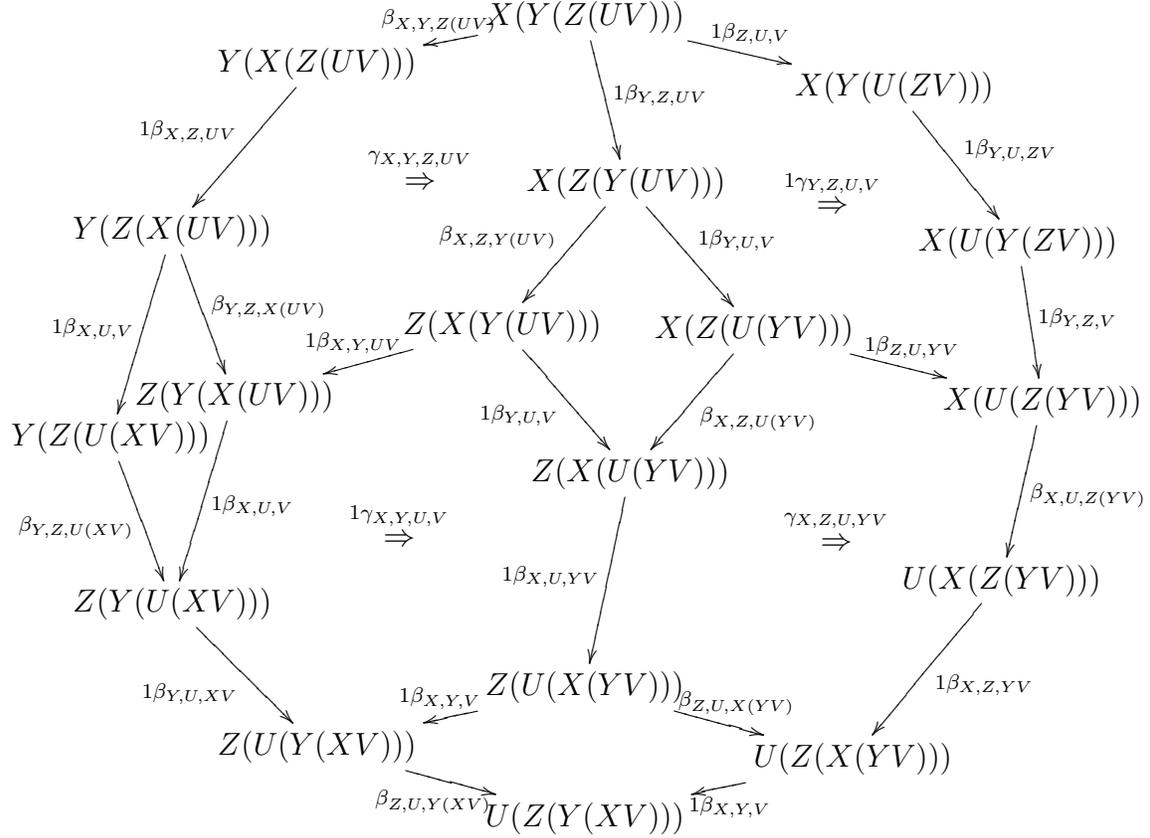
$$
\xygraph{ !{0;/r6.5pc/:;/u7.5pc/::} []*+{X(Y(Z(UV)))}
(
:[l(1.3)d(.2)]*+{Y(X(Z(UV)))} _{\beta_{X,Y,Z(UV)}}
 (
 :[l(.7)d(.7)]*+{Y(Z(X(UV)))} _{1\beta_{X,Z,UV}}
 :[l(.3)d(.8)]*+{Y(Z(U(XV)))} _{1\beta_{X,U,V}}
  (
  :[d(.9)r(.3)]*+{Z(Y(U(XV)))} _{\beta_{Y,Z,U(XV)}}
  :[d(.7)r(.7)]*+{Z(U(Y(XV)))} _{1\beta_{Y,U,XV}}
  :[d(.3)r(1.3)]*+{U(Z(Y(XV)))}="b" _{\beta_{Z,U,Y(XV)}}
  ,
  :[r(1.2)d(.3)]*+{Y(U(Z(XV))}="lm" ^{1\beta_{Z,U,XV}}
  :[r(.7)d(.65)]*+{U(Y(Z(XV))}="bm" _{\beta_{Y,U,Z(XV)}}
  :"b" ^{1\beta_{Y,Z,XV}}
  ,
  :@{}[d(1)r(1)]*{\scriptstyle{\gamma_{Y,Z,U,XV}}}
  :@{}[d(.08)]*+{\Rightarrow}
  )
 ,
 :[r(1)d(.3)]*+{Y(X(U(ZV)))}="tm" ^(.6){1\beta_{Z,U,V}}
 :[d(.85)l(.1)]*+{Y(U(X(ZV)))} ^{1\beta_{X,U,ZV}}
  (
  :"lm" _{1\beta_{X,Z,V}}
  ,
  :[r(.7)d(.66)]*+{U(Y(X(ZV)))}="rm" ^{\beta_{Y,U,X(ZV)}}
  :"bm" ^{1\beta_{X,Z,V}}
  )
 ,
 :@{}[d(.8)r(.1)]*{\scriptstyle{1\gamma_{X,Z,U,V}}}
 :@{}[d(.1)]*+{\Rightarrow}
 )
,
:[r(1.5)d(.3)]*+{X(Y(U(ZV)))} ^{1\beta_{Z,U,V}}
 (
 :"tm" _{\beta_{X,Y,U(ZV)}}
 ,
 :[d(.65)r(.6)]*+{X(U(Y(ZV)))} ^{1\beta_{Y,U,ZV}}
  (
  :[d(.9)l(.45)]*+{U(X(Y(ZV)))} _{\beta_{X,U,Y(ZV)}}
   (
   :"rm" _{1\beta_{X,Y,ZV}}
   ,
   :[d(.8)r(.4)]*+{U(X(Z(YV)))}="r" _{1\beta_{Y,Z,V}}
   :[d(.7)l(.55)]*+{U(Z(X(YV)))} _{1\beta_{X,Z,YV}}
   :"b" _{1\beta{X,Y,V}}
   ,
   :@{}[d(.8)l(.8)]*{\scriptstyle{1\gamma_{X,Y,Z,V}}}
   :@{}[d(.1)]*+{\Rightarrow}
   )
  ,
  :[d(.8)r(.4)]*+{X(U(Z(YV)))} ^{1\beta_{Y,Z,V}}
  :"r" ^{\beta_{X,U,Z(YV)}}
  )
 ,
 :@{}[d(.8)l(.8)]*{\scriptstyle{\gamma_{X,Y,U,ZV}}}
 :@{}[d(.1)]*+{\Rightarrow}
 )
)
}$$
and
$$
\xygraph{ !{0;/r6.5pc/:;/u7.5pc/::} []*+{X(Y(Z(UV)))}
(
:@{}[d(.6)l(.8)]*{\scriptstyle{\gamma_{X,Y,Z,UV}}}
:@{}[d(.1)]*+{\Rightarrow}
, 
:@{}[d(.7)r(1.2)]*{\scriptstyle{1\gamma_{Y,Z,U,V}}}
:@{}[d(.1)]*+{\Rightarrow}
 ,
:[l(1.3)d(.2)]*+{Y(X(Z(UV)))} _{\beta_{X,Y,Z(UV)}}
:[l(.7)d(.7)]*+{Y(Z(X(UV)))} _{1\beta_{X,Z,UV}}
 (
 :[l(.3)d(.87)]*+{Y(Z(U(XV)))} _{1\beta_{X,U,V}}
 :[r(.3)d(.7)]*+{Z(Y(U(XV)))}="l" _{\beta_{Y,Z,U(XV)}}
 :[d(.6)r(.7)]*+{Z(U(Y(XV)))}="lb" _{1\beta_{Y,U,XV}}
 :[d(.3)r(1.3)]*+{U(Z(Y(XV)))}="b" _{\beta_{Z,U,Y(XV)}}
 ,
 :[r(.3)d(.7)]*+{Z(Y(X(UV)))}="lm" ^{\beta_{Y,Z,X(UV)}}
 :"l" ^{1\beta_{X,U,V}}
 )
,
:[d(.7)r(.2)]*+{X(Z(Y(UV)))} ^{1\beta_{Y,Z,UV}}
 (
 :@{}[d(1.4)l(1.1)]*{\scriptstyle{1\gamma_{X,Y,U,V}}}
 :@{}[d(.1)]*+{\Rightarrow}
 ,
 :[d(.6)l(.6)]*+{Z(X(Y(UV)))} _{\beta_{X,Z,Y(UV)}}
  (
  :[d(.62)r(.62)]*+{Z(X(U(YV)))}="m" _{1\beta_{Y,U,V}}
  :[d(.9)l(.22)]*+{Z(U(X(YV)))} _{1\beta_{X,U,YV}}
   (
   :"lb" _{1\beta_{X,Y,V}}
   ,
   :[d(.3)r(1.3)]*+{U(Z(X(YV)))}="rb" ^{\beta_{Z,U,X(YV)}}
   :"b" ^{1\beta_{X,Y,V}}
   )
  ,
  :"lm" _{1\beta_{X,Y,UV}}
  )
 ,
 :[d(.62)r(.62)]*+{X(Z(U(YV)))} ^{1\beta_{Y,U,V}}
  (
  :@{}[d(.8)r(.4)]*{\scriptstyle{\gamma_{X,Z,U,YV}}}
  :@{}[d(.1)]*+{\Rightarrow}
  ,
  :"m" ^{\beta_{X,Z,U(YV)}}
  ,
  :[r(1.4)d(.3)]*+{X(U(Z(YV)))}="r" ^{1\beta_{Z,U,YV}}
  :[d(.75)l(.2)]*+{U(X(Z(YV)))} ^{\beta_{X,U,Z(YV)}}
  :"rb" ^{1\beta_{X,Z,YV}}
  )
 )
,
:[r(1.5)d(.3)]*+{X(Y(U(ZV)))} ^{1\beta_{Z,U,V}}
:[d(.65)r(.6)]*+{X(U(Y(ZV)))} ^{1\beta_{Y,U,ZV}}
:"r" ^{1\beta_{Y,Z,V}}
}$$

\bre
Note that 1-, 2-, and 3-dimensional coherences (1-morphisms $\beta$, 2-cells $\gamma$ and the above coherence for them) naturally take 
	the
shapes of the 1-, 2-, and 3-dimensional permutohedra \cite{zi}. 
\ere

A functor $\cF:\bC\to\bD$ between b-bicategories is a {\em b-functor} if it comes equipped with
a pseudo-natural in $X,Y\in\bC$ collection of equivalences
\beq\lb{FXY}
	F_{X,Y}:
F(X\ot Y)\to F(X)\ot F(Y)\ ,
\eeq
and a natural in $X,Y,Z\in\bC$ collection of 2-isomorphisms
$$\xymatrix{ 
F(X\ot (Y\ot Z))\ar[rr]^{F_{X,Y\ot Z}} \ar[dd]_{F(\beta_{X,Y,Z})} && F(X)\ot F(Y\ot Z) \ar[rr]^{1\ot F_{Y,Z}} && F(X)\ot (F(Y)\ot F(Z)) \ar[dd]^{\beta_{F(X),F(Y),F(Z)}}\\ && \Uparrow{\scriptstyle \phi_{X,Y,Z}} &&\\
F(Y\ot(X\ot Z)) \ar[rr]^{F_{Y,X\ot Z}} && F(Y)\ot F(X\ot Z) \ar[rr]^{1\ot F_{X,Z}} && F(Y)\ot(F(X)\ot F(Z))}$$
such that the pasting of 2-cells
$$
\xygraph{ !{0;/r8.4pc/:;/u7.5pc/::} []*+{F(X(Y(ZW)))}
(
:[r(1)u(.5)]*+{F(Y(X(ZW)))} ^{F(\beta_{X,Y,ZW})}
:[r(1.5)]*+{F(Y(Z(XW)))} ^{F(1\beta_{X,Z,W})}
:[r(1)d(.5)]*+{F(Z(Y(XW)))}="tr" ^{F(\beta_{Y,Z,XW})}
:[d]*+{F(Z)F(Y(XW))}="mr" ^{F_{Z,Y(XW)}}
:[d]*+{F(Z)(F(Y)F(XW))} ^{1F_{Y,XW}}
:[d]*+{F(Z)(F(Y)(F(X)F(W)))}="br" ^{11F_{X,W}}
,
:[r(1)d(.5)]*+{F(X(Z(YW)))} ^{F(1\beta_{Y,Z,W})}
 (
 :[r(1.5)]*+{F(Z(X(YW)))} ^{F(\beta_{X,Z,YW})}
  (
  :"tr" ^{F(1\beta_{X,Y,W})}
  ,
  :[d]*+{F(Z)F(X(YW))} _{F_{Z,X(YW)}}
   (
   :"mr" ^{1F(\beta_{X,Y,W})}
   ,
   :[d]*+{F(Z)(F(X)F(YW))}="rm" _{1F_{X,YW}}
   :[d]*+{F(Z)(F(X)(F(Y)F(W)))}="rb" _{11F_{Y,W}}
   :"br" _(.6){\quad 1\beta_{F(X),F(Y),F(W)}}
   :@{}[u(.8)l(.5)]*+{\Leftarrow}
   :@{}[d(.1)]*+{\scriptstyle{1\phi_{X,Y,W}}}
   )
  )
 ,
 :[d]*+{F(X)F(Z(YW))}="lm" ^{F_{X,Z(YW)}}
 :[d]*+{F(X)(F(Z)F(YW))} ^{1F_{Z,YW}}
  (
  :"rm" _{\beta_{F(X),F(Z),F(YW)}}
  :@{}[u(1)l(.7)]*+{\Leftarrow}
  :@{}[d(.1)]*+{\scriptstyle{\phi_{X,Z,YW}}} 
  ,
  :[d]*+{F(X)(F(Z)(F(Y)F(W)))}="b" ^{11F_{Y,W}}
  :"rb" _{\beta_{F(X),F(Z),F(Y)F(W)}}
  )
 )
,
:[d]*+{F(X)F(Y(ZW))} _{F_{X,Y(ZW)}}
 (
 :"lm" ^{1F(\beta_{Y,Z,W})}
 ,
 :[d]*+{F(X)(F(Y)F(ZW))} _{1F_{Y,ZW}}
 :[d]*+{F(X)(F(Y)(F(Z)F(W)))} _{11F_{Z,W}}
 :"b" _{1\beta_{F(Y),F(Z),F(W)}}
 :@{}[u(1.3)l(.5)]*+{\Leftarrow}
 :@{}[d(.1)]*+{\scriptstyle{1\phi_{Y,Z,W}}}
 )
,
:@{}[r(1.8)]*{\Downarrow\scriptstyle{F(\gamma_{X,Y,Z,W})}}
)
}$$
coincides with the pasting
$$
\xygraph{ !{0;/r8.4pc/:;/u7.5pc/::} []*+{F(X(Y(ZW)))}
(
:[r(1)u(.5)]*+{F(Y(X(ZW)))} ^{F(\beta_{X,Y,ZW})}
 (
 :[r(1.5)]*+{F(Y(Z(XW)))} ^{F(1\beta_{X,Z,W})}
  (
  :[r(1)d(.5)]*+{F(Z(Y(XW)))} ^{F(\beta_{Y,Z,XW})}
  :[d]*+{F(Z)F(Y(XW))} ^{F_{Z,Y(XW)}}
  :[d]*+{F(Z)(F(Y)F(XW))}="rm" ^{1F_{Y,XW}}
  :[d]*+{F(Z)(F(Y)(F(X)F(W)))}="br" ^{11F_{X,W}}
  ,
  :[d]*+{F(Y)F(Z(XW))}="tr" _{F_{Y,Z(XW)}}
  :[d]*+{F(Y)(F(Z)F(XW))} _{1F_{Z,XW}}
   (
   :"rm" ^{\beta_{F(Y),F(Z),F(XW)}}
   :@{}[u(1.3)l(.5)]*+{\Leftarrow}
   :@{}[d(.1)]*+{\scriptstyle{\phi_{Y,Z,XW}}}
   ,
   :[d]*+{F(Y)(F(Z)(F(X)FW)))}="mr" _{11F_{X,W}}
   :"br" _{\beta_{F(Y),F(Z),F(X)F(W)}}
   )
  )
 ,
 :[d]*+{F(Y)F(X(ZW))} ^{F_{Y,X(ZW)}}
  (
  :"tr" ^{1F(\beta_{X,Z,W})}
  ,
  :[d]*+{F(Y)F(X)F(ZW))}="m" ^{1F_{X,ZW}}
  :[d]*+{F(Y)(F(X)(F(Z)F(W)))}="ml" ^{11F_{Z,W}}
  :"mr" _{1\beta_{F(X),F(Z),F(W)}}
  :@{}[u(1)l(.7)]*+{\Leftarrow}
  :@{}[d(.1)]*+{\scriptstyle{1\phi_{X,Z,W}}}
  )
 )
,
:[d]*+{F(X)F(Y(ZW))} _{F_{X,Y(ZW)}}
:[d]*+{F(X)(F(Y)F(ZW))} _{1F_{Y,ZW}}
 (
 :[d]*+{F(X)(F(Y)(F(Z)F(W)))} _{11F_{Z,W}} 
  (
  :[d(.5)r]*+{F(X)(F(Z)(F(Y)F(W)))} _{1\beta_{F(Y),F(Z),F(W)}}
  :[r(1.5)]*+{F(X)(F(Z)(F(Y)F(W)))} _{\beta_{F(X),F(Z),F(Y)F(W)}}
  :"br" _(.6){\quad 1\beta_{F(X),F(Y),F(W)}}
  :@{}[l(1.7)]*+{\Downarrow\scriptstyle{\gamma_{F(X),F(Y),F(Z),F(W)}}} 
  ,
  :"ml" _{\quad\beta_{F(X),F(Y),F(Z)F(W)}}
  )
,
:"m" ^{\beta_{F(X),F(Y),F(ZW)}}
:@{}[u(.7)l(.5)]*+{\Leftarrow} 
:@{}[d(.1)]*+{\scriptstyle{\phi_{X,Y,ZW}}}
 )
)
}$$

The composition of b-functors $\xymatrix{\bC\ar[r]^F & \bD \ar[r]^G & \bE}$ is again a b-functor with   isomorphisms $(G\circ F)_{X,Y}:(G\circ F)(X\ot Y)\to (G\circ F)(X)\ot (G\circ F)(Y)$ defined as in \eqref{cobf}. The natural 2-cells for the composition $G\circ F$ are defined by the pasting
\beq\lb{c2c}
\xygraph{ !{0;/r7.5pc/:;/u5.5pc/::} []*+{(GF)(X(YZ))}
(
:[rr]*+{GF(X)GF(YZ)}="tm" ^{(GF)_{X,YZ}}
 (
 :[rr]*+{GF(X)(GF(Y)GF(Z))}="tr" ^{1(GF)_{Y,Z}}
 :[d(5)]*+{GF(Y)(GF(X)GF(Z))}="br" ^{\beta_{GF(X),GF(Y),GF(Z)}}
 ,
 :[rd]*+{(GF)(X)G(F(Y)F(Z))}="mr" ^{1G(F_{Y,Z})}
 :"tr" _{1G_{F(Y),F(Z)}}
 )
,
:[rd]*+{G(F(X)F(YZ))} _{G(F_{X,YZ})}
 (
 :"tm" ^{G_{F(X),F(YZ)}}
 ,
 :[rd]*+{G(F(X)(F(Y)F(Z)))} _{G(1F_{Y,Z})}
  (
  :"mr" _{G_{F(X),F(Y)F(Z)}}
  ,
  :[d]*+{G(F(Y)(F(X)F(Z)))}="m" _(.4){G(\beta_{F(X),F(Y),F(Z)})}
  :[rd]*+{(GF)(Y)G(F(X)F(Z))}="mb" ^{G_{F(Y),F(X)F(Z)}}
  :"br" ^{1G_{F(X),F(Z)}}
  :@{}[u(2.3)l(.8)]*+{\Uparrow\scriptstyle{\psi_{F(X),F(Y),F(Z)}}}
  )
 )
,
:[d(5)]*+{(GF)(Y(XZ))} _{(GF)(\beta_{X,Y,Z})}
 (
 :@{}[u(2.3)r(.8)]*+{\Uparrow\scriptstyle{G(\phi_{X,Y,Z})}}
 ,
 :[rr]*+{(GF)(Y)(GF)(XZ)}="b" _{(GF)_{Y,XZ}}
  (
  :"br" _{1(GF)_{X,Z}}
  ,
  :"mb" _{1G(F_{X,Z})}
  )
 ,
 :[ru]*+{G(F(Y)F(XZ))} ^{G(F_{Y,XZ})}
  (
  :"b" ^{G_{F(Y),F(XZ)}}
  ,
  :"m" _{G(1F_{X,Z})}
  )
 )
)
}
\eeq
Here $\phi$ is the structure 2-cell of $F$ and $\psi$ is the structure 2-cell of $G$.
It is straightforward to verify that this 2-isomorphism satisfies the coherence axiom. 

A pseudo-natural transformation $a:F\to G$ of two b-functors $F,G:\bC\to\bD$ between b-bicategories is {\em b-transformation} if it comes equipped with a collection of 2-isomorphisms 
$$\xygraph{ !{0;/r4.5pc/:;/u3.5pc/::} []*+{F(X\ot Y)}
(
:[dd]*+{G(X\ot Y)} _{a_{X\ot Y}}
:[rr]*+{G(X)\ot G(Y)}="e" _{G_{X,Y}}
,
:[rr]*+{F(X)\ot F(Y)} ^{F_{X,Y}}
:"e" ^{a_X\ot a_Y}
,
:@{}[r(1.1)d]*+{\Downarrow\scriptstyle{\alpha_{X,Y}}}
)
}$$
natural in $X,Y\in\bC$ such that the pasting
\beq\lb{pntcoh1}
\xygraph{ !{0;/r5.5pc/:;/u5.5pc/::} []*+{F(X(YZ))}
(
:[dl]*+{F(Y(XZ))} _{F(\beta_{X,Y,Z})}
 (
 :[dr]*+{G(Y(XZ))} _{a_{Y(XZ)}}
 :[rr]*+{G(Y)G(XZ)}="b" _{G_{Y,XZ}}
 :[rr]*+{G(Y)(G(X)G(Z))}="br" _{1G_{X,Z}}
 ,
 :[rr]*+{F(Y)F(XZ)} ^{F_{Y,XZ}}
  (
  :"b" ^{a_Ya_{XZ}}
  ,
  :[rr]*+{F(Y)(F(X)F(Z))}="m" ^{1F_{X,Z}}
  :"br"  ^{a_Ya_Xa_Z}
  ,
  :@{}[r(1.5)d(.5)]*+{\Downarrow\scriptstyle{1\alpha_{X,Z}}}
  )
 ,
 :@{}[r(1.5)d(.5)]*+{\Downarrow\scriptstyle{\alpha_{Y,XZ}}}
 )
,
:[rr]*+{F(X)F(YZ)} ^{F_{X,YZ}}
:[rr]*+{F(X)(F(Y)F(Z))} ^{1F_{Y,Z}}
 (
 :"m" _{\beta_{F(X),F(Y),F(Z)}}
 ,
 :[rd]*+{G(X)(G(Y)G(Z))} ^{a_Xa_Ya_Z}
 :"br" ^{\beta_{G(X),G(Y),G(Z)}}
 )
,
:@{}[r(1.6)d(.5)]*+{\Uparrow\scriptstyle{\phi_{X,Y,Z}}}
)
}\eeq
coincides with
\beq\lb{pntcoh2}
\xygraph{ !{0;/r5.5pc/:;/u5.5pc/::} []*+{F(X(YZ))}
(
:[dl]*+{F(Y(XZ))} _{F(\beta_{X,Y,Z})}
:[dr]*+{G(Y(XZ))}="bl" _{a_{Y(XZ)}}
:[rr]*+{G(Y)G(XZ)} _{G_{Y,XZ}}
:[rr]*+{G(Y)(G(X)G(Z))}="br" _{1G_{X,Z}}
,
:[dr]*+{G(X(YZ))} _{a_{X(YZ)}}
 (
 :"bl" _{G(\beta_{X,Y,Z})}
 ,
 :[rr]*+{G(X)G(YZ)}="m" ^{G_{X,YZ}}
 :[rr]*+{G(X)(G(Y)G(Z))}="mr" ^{1G_{Y,Z}}
 :"br" ^{\beta_{G(X),G(Y),G(Z)}}
 ,
 :@{}[r(1.8)d(.5)]*+{\Uparrow\scriptstyle{\psi_{X,Y,Z}}}
 )
,
:[rr]*+{F(X)F(YZ)} ^{F_{X,YZ}}
 (
 :"m" ^{a_Xa_{YZ}}
 ,
 :[rr]*+{F(X)(F(Y)F(Z))} ^{1F_{Y,Z}}
 :"mr" ^{a_Xa_Ya_Z}
 ,
 :@{}[r(1.5)d(.5)]*+{\Downarrow\scriptstyle{1\alpha_{Y,Z}}}
 )
,
:@{}[r(1.5)d(.5)]*+{\Downarrow\scriptstyle{\alpha_{X,YZ}}}
)
 }\eeq
for any $X,Y,Z\in\bC$. 

A modification 
$$\xymatrix{F(X) \ar@/^20pt/[rr]^{a_X} \ar@/_20pt/[rr]_{a'_X} & \Downarrow\scriptstyle{m_X} & G(X)}$$
between pseudo-natural b-transformations is a {\em b-modification} if the pasting
$$\xygraph{ !{0;/r4.5pc/:;/u3.5pc/::} []*+{F(X\ot Y)}
(
:@/^20pt/[dd]*+{G(X\ot Y)}="b" ^{a_{X\ot Y}}
:[rr]*+{G(X)\ot G(Y)}="e" _{G_{X,Y}}
,
:[rr]*+{F(X)\ot F(Y)} ^{F_{X,Y}}
:@/^20pt/"e" ^{a_X\ot a_Y}
,
:@{}[r(1.3)d(.9)]*+{\Leftarrow}
:@{}[d(.2)]*+{\scriptstyle{\alpha_{X,Y}}}
,
:@/_20pt/"b" _{a'_{X\ot Y}}
,
:@{}[d(.9)]*+{\Leftarrow}
:@{}[d(.2)]*+{\scriptstyle{m_{X\ot Y}}}
)
}$$
coincides with
$$\xygraph{ !{0;/r4.5pc/:;/u3.5pc/::} []*+{F(X\ot Y)}
(
:@/_20pt/[dd]*+{G(X\ot Y)}="b" _{a'_{X\ot Y}}
:[rr]*+{G(X)\ot G(Y)}="e" _{G_{X,Y}}
,
:[rr]*+{F(X)\ot F(Y)} ^{F_{X,Y}}
 ( 
 :@/^20pt/"e" ^{a_X\ot a_Y}
 ,
 :@/_20pt/"e" _{a'_X\ot a'_Y}
 ,
 :@{}[d(.9)]*+{\Leftarrow}
 :@{}[d(.2)]*+{\scriptstyle{m_X\ot m_Y}} 
 )
,
:@{}[r(.6)d(.9)]*+{\Leftarrow}
:@{}[d(.2)]*+{\scriptstyle{\alpha'_{X,Y}}}
)
}$$
for any $X,Y\in\bC$. 

\bre\lb{tbbc}
It is straightforward to see that pseudo-natural b-transformations and b-modifications are composable and are compatible with the composition of b-functors.
In other words b-bicategories, b-functors, pseudo-natural b-transformations and b-modifications form a tri-category \cite{gps}.
In particular the 2-category of b-endofunctors, pseudo-natural b-transformations and b-modifications for a fixed b-bicategory
	is
a monoidal 2-category.
\ere

\subsection{One object b-bicategories and Zamolodchikov's tetrahedron equation}

Let $B$ be an object of a monoidal category $\E$ equipped with a half braiding $t_{B,X} : B \ot X \to X \ot B$ (i.e.\ $B$ is in the monoidal centre of $\E$), such that $t$ is symmetric in the sense that $(t_{B,B})^2 = 1$.
We say that an endomorphism $Z:B^{\ot 3}\to B^{\ot 3}$ satisfies the {\em Zamolodchikov tetrahedron equation} (see e.g. \cite[section 1.7]{kv} or \cite{bs}) if the following holds in $\E(B^{\ot 6},B^{\ot 6})$:
$$Z_{124}Z_{135}Z_{236}Z_{456} = Z_{456}Z_{236}Z_{135}Z_{124}\ .$$
Here 
$Z_{124}=t^{-1}_3(Z\ot 1)t_3$, etc.

\medskip

By a {\em weak} b-bicategory we will mean a b-bicategory with $\beta$ being not necessarily equivalences but just  1-morphisms.

Here we look at a weak b-bicategory $\bC$ with only one object $I$. By $\E=\bC(I,I)$ we denote the 
category of endomorphisms, i.e.\ the category with objects being morphisms $I\to I$ in $\bC$ and with morphisms being 2-cells in $\bC$. Composition of morphisms in $\bC$ makes $\E$ a monoidal category.
We express the composition of 1-morphisms $I \xrightarrow{A} I \xrightarrow{B} I \xrightarrow{C} I$ in $\bC(I,I)$ as tensor product $A \ot B \ot C$ in $\E$ (rather than in the opposite order). This will be important when translating pasting diagrams into compositions of morphisms in $\E$.
Since $I\ot(I\ot I)=I$ the 1-morphism $\beta_{I,I,I}$ is an object $B\in\E$. 
Pseudo-naturality 2-cells \eqref{nat} for $\beta$ (with only one of the 1-morphisms $f,g,h$ being non-identical) correspond to a half-braiding $t_{B,X}:B\ot X\to X\ot B$,
	which in the following we will assume to be symmetric in the above sense.

\bth\lb{oobb}
The data of a weak b-bicategory $\bC$ with only one object $I$ and with category of endomorphisms $\E=\bC(I,I)$ 
	such that $t_{B,B}^2 = 1$
is equivalent to a pair $(B,Z)$ consisting of an object $B\in\E$ and a solution $Z\in\E(B^{\ot 3},B^{\ot 3})$ of the Zamolodchikov tetrahedron equation.
\eth
\bpf 
Upon the identification $I\ot(I\ot(I\ot I))=I$ both $1\ot\beta_{I,I,I}$ and $\beta_{I,I,I\ot I}$ coincide with $B$. Thus the source and the target of $\gamma_{I,I,I,I}$ are $B^{\ot 3}$ and $\gamma_{I,I,I,I}$ is an isomorphism 
$S:B^{\ot 3}\to B^{\ot 3}$ in the category $\E$. 
\nl
Each of the pasting diagrams of the coherence condition for $\gamma$ gives a composition of tensor products of the identity with $S$ (possibly conjugated with braiding). Since all the paths (the sequences of successive 1-morphisms) of each pasting diagram have length six each factor in the products is an automorphism of $B^{\ot 6}$. Note that squares of the form
$$\xygraph{ !{0;/r8.5pc/:;/u6.5pc/::} []*+{I^{\ot 6}}
(
:[r]*+{I^{\ot 6}} ^{\beta_{I,I,I^{\ot 3}}}
:[d]*+{I^{\ot 6}}="e" ^{1\ot\beta_{I,I,I}}
,
:[d]*+{I^{\ot 6}} _{1\ot\beta_{I,I,I}}
:"e" _{\beta_{I,I,I^{\ot 3}}}
,
:@{}[r(.6)d(.5)]*+{\Downarrow \scriptstyle{b_{1,1,\beta_{I,I,I}}}}
)
}$$ correspond to the half-braiding $t_{B,B}:B\ot B\to B\ot B$. 
Thus the equational form of two pasting diagrams (specialised at $X=Y=Z=U=V=I$) is as follows
\beq\lb{gamcoh}
t_3S_{456}S_{234}(t_1t_4)S_{234}S_{456} = S_{123}S_{345}(t_2t_5)S_{345}S_{123}t_3\ .
\eeq
Substituting $S = t_1t_2t_1Z$ we get (by pulling all $t$'s to the right and changing the indices of the $Z$'s accordingly)
\begin{align*}
& t_3t_4t_5t_4t_2t_3t_2t_1t_4t_2t_3t_2t_4t_5t_4 \, Z_{124}Z_{135}Z_{236}Z_{456}  \\
& = t_1t_2t_1t_3t_4t_3t_2t_5t_3t_4t_3t_1t_2t_1t_3 \, Z_{456}Z_{236}Z_{135}Z_{124}\ .
\end{align*}
Finally since 
	the words in $t$'s on the two sides are equal
in the symmetric group $S_6$, we get the Zamolodchikov tetrahedron equation. 
\epf

By a {\em weak b-functor} we will mean a b-functor $F:\bC\to\bD$ between (weak) b-bicategories with $F_{X,Y}$ being not necessarily equivalences but just 1-morphisms.

Now we look at weak b-endofunctors of a weak b-bicategory $\bC$ with only one object $I$. 
The 1-morphism $F_{I,I}$ in \eqref{FXY} with $F=Id$ is an object $C\in\E$. 
Pseudo-naturality 2-cells  for $F_{I,I}$  gives a half-braiding $t_{C,X}:C\ot X\to X\ot C$. 
In what follows we assume that the half-braiding for $B$, $t_{B,C}:B\ot C\to C\ot B$, and the half-braiding for $C$, $t_{C,B}:C\ot B\to B\ot C$, are inverse to each other. 

\bpr\lb{ibf}
Let $\bC$, $I$, $\E$, as well as $B$ and $Z$ be as in Theorem \ref{oobb}.
The structure of weak b-functor on the identity functor $Id:\bC\to\bC$  
	such that the half-braidings $t_{B,C}$ and $t_{C,B}$ are each other's inverses
amounts to a pair $(C,L)$ consisting of an object $C\in\E$ and an automorphism $L:C^{\ot 2}\ot B\to C^{\ot 2}\ot B$ such that
\beq\lb{lze}L_{124}L_{135}L_{236}Z_{456} = Z_{456}L_{236}L_{135}L_{124}\eeq
in $\E(C^{\ot 3}\ot B^{\ot 3},C^{\ot 3}\ot B^{\ot 3})$. 
\epr
\bpf
All vertices in the diagram for $\phi_{I,I,I}$ are equal to $I$ and
both $1\ot F_{I,I}$ and $F_{I,I\ot I}$ coincide with $C$.
Replacing also $\beta_{I,I,I}=B \in \E$, we see that
the source of $\phi_{I,I,I}$ is $B\ot C^{\ot 2}$ 
while the target is $C^{\ot 2}\ot B$. 
The 2-cell $\phi_{I,I,I}$ is an isomorphism $M:B\ot C^{\ot 2}\to C^{\ot 2}\ot B$ in the category $\E$. 

In the pasting diagrams of the coherence condition for $\phi_{I,I,I}$ the cells $F(\gamma_{I,I,I,I})$ and $\gamma_{F(I),F(I),F(I),F(I)}$ are given by the morphism $S:B^{\ot 3}\to B^{\ot 3}$ 
	from the proof of Theorem \ref{oobb}.
Cells labelled by $\phi$ (and their tensor products with the identity) correspond to instances of the isomorphism $M:B\ot C^{\ot 2}\to C^{\ot 2}\ot B$. Squares with opposite sides labelled by $\beta$ and $F$ correspond 
	to one of the half-braidings $t_{B,C}:B\ot C\to C\ot B$ or $t_{C,B}:C\ot B\to B\ot C$. 
The coherence for $\phi_{I,I,I}$ can be reformulated as the commutativity of the diagram
$$\xymatrix{B^{\ot 3}C^{\ot 3} \ar[r]^{S_{123}} \ar[d]_{M_{345}} & B^{\ot 3}C^{\ot 3} \ar[r]^{t_3} & BBCBCC \ar[r]^{M_{456}} & BBC^{\ot 3}B \ar[r]^{M_{234}} & BCCBCB \ar[d]^{t_1t_4}  \\ 
BBCCBC \ar[d]_{t_2t_5} &&&& CBCCBB \ar[d]^{M_{234}} \\
BCBCCB \ar[r]^{M_{345}} & BC^{\ot 3}BB \ar[r]^{M_{123}} & CCBCBB \ar[r]^{t_3} & C^{\ot 3}B^{\ot 3} \ar[r]^{S_{456}} & C^{\ot 3}B^{\ot 3}}$$
or, in equational form,
$$M_{234}(t_1t_4)M_{234}M_{456}t_3S_{123} = S_{456}t_3M_{123}M_{345}(t_2t_5)M_{345}\ .$$
Define $L =  t_1t_2t_1M^{-1}$. Rewriting the above equation as
$$t_3M^{-1}_{456}M^{-1}_{234}(t_1t_4)M^{-1}_{234}S_{456} = S_{123}M^{-1}_{345}(t_2t_5)M^{-1}_{345}M^{-1}_{123}t_3$$ 
one sees that this has the same index structure as \eqref{gamcoh} in the proof of Theorem \ref{oobb}, so that by the same argument as used there, one arrives at \eqref{lze}.
\epf

\bre
Equation \eqref{lze} appears in integrable three-dimensional lattice models of statistical mechanics \cite{bs}. 
Our interpretation allows one to see a 2-categorical structure on solutions of the equation \eqref{lze}. 
Indeed, b-functors from Proposition \ref{ibf}\ are composable and form a monoidal 2-category (the sub-2-category from Remark \ref{tbbc}).
The corresponding structure of a monoidal bicategory can be described in terms of pairs $(C,L)$:
\begin{itemize}\setlength{\leftskip}{-1em}
\item[-]
A morphism $(C,L)\to (C',L')$ is an object $D\in\E$ together with an isomorphism 
	$d:C\ot D^{\ot 2} \to D\ot C'$ satisfying
	coherence conditions coming from \eqref{pntcoh1} and \eqref{pntcoh2}.
\item[-]
A 2-morphism $(D,d)\to (D',d')$ (between two morphisms $(C,L)\to (C',L')$) is a morphism $f:D\to D'$ in $\E$ such that
the following diagram commutes:
$$
\xymatrix{C\ot D^{\ot 2} \ar[r]^d \ar[d]_{1\ot f^{\ot 2}} & D \ot C' \ar[d]^{f \ot 1} \\ C\ot D'^{\ot 2}\ar[r]^{d'}  & D'\ot C'}
$$
\item[-]
The composition of b-functors corresponds to the tensor product $(C,L)\ot (C',L')$ defined as $(C\ot C',L|L')$, where $L|L'$ is the composition
 $t_2 L_{125} L'_{345}t_2$.
In particular, $L|L'$ is a solution of the equation \eqref{lze}, which is not obvious.
\end{itemize}
\ere

\section{b-cohomology} \label{sec:b-cohomology}

The coherence conditions of b-functors, b-categories and their bicategorical versions motivate the definition of a certain cochain complex, whose cohomology we call b-cohomology.

\subsection{Definition}\lb{bcoh}

Let $A$ be a b-magma and $B$ be an abelian group (with multiplicatively and additively written operations respectfully).
Let $C^n(A,B) = Maps(A^{\times n},B)$ be the (additive) abelian group of maps $A^{\times n}\to B$.
For $n\geq 2$ define the map 
$$d:C^n(A,B)\to C^{n+1}(A,B)$$
by 
$$d(c)(x_1,...,x_n,x) = \sum_{i=1}^n(-1)^{i}\delta_i(c)(x_1,...,x_n,x) \ ,$$
where  
$\delta_i(c)(x_1,...,x_n,x) = c(x_1,...,\widehat{x_i},...,x_{n},x_ix) - c(x_1,...,\widehat{x_i},...,x_{n},x)$
and $\hat{\ }$ above a variable means its omission.

\medskip

For example, for $q\in C^2(A,B)$, $r \in C^3(A,B)$ and  $s \in C^4(A,B)$,
\begin{align*}
d(q)(x,y,z) &= - q(y,xz) + q(y,z) +  q(x,yz) - q(x,z)   \ ,\\
d(r)(x,y,z,w) &= - r(y,z,xw) + r(y,z,w) + r(x,z,yw) - r(x,z,w) \\ & ~~~\, - r(x,y,zw) + r(x,y,w)  \ ,\\
d(s)(x,y,z,u,v) &= 	- s(y,z,w,xv) + s(y,z,u,v)
				+ s(x,z,u,yv) - s(x,z,u,v) \\ & ~~~\,
				- s(x,y,u,zv) + s(x,y,u,v)
				+ s(x,y,z,uv) - s(x,y,z,v) \ .
\end{align*}

\ble
The map $d:C^n(A,B)\to C^{n+1}(A,B)$ is a differential, i.e.\ $d^2=0$.
\ele
\bpf
First we need to show that $\delta_i\circ \delta_j = \delta_{j+1}\circ \delta_i$ for $i\leq j$. Indeed,
\begin{align*}
&\delta_i(\delta_j(c))(x_1,...,x_{n+1},x) \\
&\quad =~ \delta_j(c)(x_1,...,\widehat{x_i},...,x_{n+1},x_ix) ~-~ \delta_j(c)(x_1,...,\widehat{x_i},...,x_{n+1},x) \\
&\quad =~ c(x_1,...,\widehat{x_i},...,\widehat{x_{j+1}},...,x_{n+1},x_{j+1}(x_ix)) ~-~ c(x_1,...,\widehat{x_i},...,\widehat{x_{j+1}},...,x_{n+1},x_ix)  \\
&\qquad ~  - c(x_1,...,\widehat{x_i},...,\widehat{x_{j+1}},...,x_{n+1},x_{j+1}x) ~+~ c(x_1,...,\widehat{x_i},...,\widehat{x_{j+1}},...,x_{n+1},x)
\end{align*}
and
\begin{align*}
&\delta_{j+1}(\delta_i(c))(x_1,...,x_{n+1},x) \\
& \quad = ~\delta_i(c)(x_1,...,\widehat{x_{j+1}},...,x_{n+1},x_{j+1}x) ~-~ \delta_i(c)(x_1,...,\widehat{x_{j+1}},...,x_{n+1},x) \\
&\quad = ~c(x_1,...,\widehat{x_i},...,\widehat{x_{j+1}},...,x_{n+1},x_i(x_{j+1}x)) ~-~ c(x_1,...,\widehat{x_i},...,\widehat{x_{j+1}},...,x_{n+1},x_{j+1}x) \\
& \qquad ~ - c(x_1,...,\widehat{x_i},...,\widehat{x_{j+1}},...,x_{n+1},x_ix) ~+~ c(x_1,...,\widehat{x_i},...,\widehat{x_{j+1}},...,x_{n+1},x) \ .
\end{align*}
The two sides coincide because of the identity $x_{j+1}(x_ix) = x_i(x_{j+1}x)$.
Now 
\begin{align*}
d^2 &= \sum_{i=1}^{n+1}\sum_{j=1}^n(-1)^{i+j}\delta_i\circ\delta_j\ = \sum_{1\leq i\leq j\leq n}(-1)^{i+j}\delta_i\circ\delta_j\ + \sum_{1\leq j<i\leq n+1}(-1)^{i+j}\delta_i\circ\delta_j\ \\
&=\sum_{1\leq i\leq j\leq n}(-1)^{i+j}\delta_{j+1}\circ\delta_i\ + \sum_{1\leq j<i\leq n+1}(-1)^{i+j}\delta_i\circ\delta_j\ \\
&= - \sum_{1\leq j<i\leq n+1}(-1)^{i+j}\delta_i\circ\delta_j\ + \sum_{1\leq j<i\leq n+1}(-1)^{i+j}\delta_i\circ\delta_j\ = 0\ .
\end{align*}
\epf

We call the cohomology $H^*_b(A,B)$ of the complex $C^*_b(A,B)=(C^*(A,B),d)$ the {\em b-cohomology} of $A$ with coefficients in $B$.

\bre
The formula defining the differential $d$ on $C^n$ makes sense for $n=1$ as well (and would read $d(a)(x,y) = -a(xy)+a(y)$), but the resulting degree 2 cohomology does not match the examples of b-magma algebras or b-functors discussed below. These examples instead suggest to extend the cochain complex $C^*_b(A,B)$ to the first degree by defining $d:C^1(A,B)\to C^2(A,B)$ as
\beq\lb{fd}d(p)(x,y) = p(x) - p(xy) + p(y) \eeq
for all $p \in C^1(A,B)$.
One verifies that indeed $d^2=0$. 
One can define $H^1(A,B) = Z^1(A,B)$.
\ere

\subsection{Realisations of low dimensional b-cohomology}\lb{blow}

Let $A$ be a b-magma. 
By a b-{\em magma algebra} over a field $k$ we mean the vector space $kA$ equipped with the bilinear multiplication defined on basis elements as $(x,y) \mapsto xy$. By construction, a b-magma algebra is a b-algebra. 

Similarly to the relation between twistings of the multiplication of group algebras and the second group cohomology  there is a relation between twistings of the multiplication of b-magma algebras and the second b-cohomology.
Namely, one can  try to twist the multiplication by a 2-cochain $q \in C^2(A,k^\times)$ and ask if $(x,y) \mapsto q(x,y) xy$ is again a b-algebra. One finds the condition
$q(y,z) q(x,yz) = q(x,z) q(y,xz)$, which amounts to $d(q)=0$.

Isomorphisms which preserve the $A$-grading act by rescaling the basis elements. Thus, the b-magma algebras obtained from $kA$ by twisting with $q$ and $q'$ are $A$-grading preserving isomorphic iff $q'(x,y)p(xy)^{-1} = q(x,y) p(x)p(y)$ for some $p \in C^1(A,k^\times)$. 
This amounts to $q' = qdp$.

Thus $H_b^2(A,k^\times)$ describes twisted b-magma algebras up to grading-preserving isomorphisms.

Similarly one can twist the identity homomorphism of b-magma algebra $kA$ by a 1-cochain $p \in C^1(A,k^\times)$ and ask if $x \mapsto p(x) x$ is again a homomorphism. One finds the condition
$p(xy) = p(x)p(y)$, which amounts to $d(p)=0$.

Thus $H_b^1(A,k^\times)$ describes grading-preserving automorphisms of the b-magma algebra $kA$.

\bigskip

One categorical level up, in analogy to the relation between braided categorical 
groups (certain braided monoidal groupoids) and abelian group cohomology \cite{js} there is a relation between the categorical b-magmas discussed in Section \ref{pbc} and b-cohomology.
Indeed, the calculations in Section \ref{pbc} can be summarised as follows.

\bpr\lb{catbmag}
Let $S$ be a b-magma. 
\nl
(i) Equivalence classes of categorical b-magmas with the set $\pi_0(\C)=S$ of isomorphisms classes of objects and typical automorphism group $\pi_1(\C)$ correspond to orbits of the automorphism group of $S$ on the b-cohomology group $H^3_b(\pi_0(\C),\pi_1(\C))$. 
\nl
(ii) The isomorphism class of a b-functor $F:\C\to\D$ between two categorical b-magmas corresponds to a homomorphism of b-magmas $\pi_0(F):\pi_0(\C)\to \pi_0(\D)$ together with an element of the b-cohomology group $H^2_b(\pi_0(\C),\pi_1(\D))$.
\epr

\bre\lb{coho}
The relation between b-categories and braided monoidal categories provides examples of homomorphisms from abelian group cohomology to b-cohomology.
Namely, for abelian groups $A$ and $B$ the assignments
\begin{align*}
C^1_{ab}(A,B)&\to C^1(A,B) \ ,~~ f\mapsto f ~ , \qquad\qquad
C^2_{ab}(A,B)\to C^2(A,B)\ ,~~ g\mapsto g \ , \\
C^3_{ab}(A,B)&\to C^3(A,B)\ ,~~ (a,c) \mapsto b \ ,
\end{align*}
where $b(x,y,z) = a(x,y,z) + c(x,y) - a(y,x,z)$, define
homomorphisms between cohomology groups:
$$H^1_{ab}(A,B)\to H_b^1(A,B) ,\qquad H^2_{ab}(A,B)\to H_b^2(A,B) ,\qquad H^3_{ab}(A,B)\to H_b^3(A,B)\ .$$
\ere

\bigskip

Finally, on the b-bicategorical level, consider the following instance of a bicategorical b-magma: the set of objects is a b-magma $S$, the set of 1-morphisms $x \to y$ ($x,y \in S$) is empty unless $x=y$. For $x=y$ there is a unique 1-morphism $\one_x$ and the 2-morphisms $\one_x \Rightarrow \one_x$ are given by a fixed abelian group $B$. Since the endomorphism categories are one-object braided monoidal categories, the braiding is necessarily trivial, i.e.\ the structure morphisms in \eqref{sli} are all identities.

We necessarily have $\beta_{x,y,z} = \one_{xyz}$, which implies that the corresponding 2-cells \eqref{nat} are identities as well. 
The structure maps $\gamma_{w,x,y,z}$ are given by a cochain $s(w,x,y,z)$ from $C^4(S,B)$.
The coherence condition for $\gamma$ becomes, for $x,y,z,u,v \in S$,
\begin{align*}
	&s(x,y,z,v) s(x,y,u,zv) s(y,z,u,xv) s(x,z,u,v)  \\
	&=
	s(y,z,u,v) s(x,z,u,yv) s(x,y,u,v) s(x,y,z,uv) \ ,
\end{align*}
in other words, $ds = 0$. Let us denote this b-bicategory by $\mathcal{B}(S,B,s)$. 

We can now ask if $\mathcal{B}(S,B,s)$ and $\mathcal{B}(S,B,s')$ are equivalent via a b-functor whose underlying functor of bicategories is the identity functor. To equip the identity functor with a b-structure we need to choose $F_{x,y}$ and $\phi_{x,y,z}$. For $F_{x,y}$ the only choice available is $\one_{xy}$, while the $\phi$'s are given by a 3-cochain $r$ subject to the coherence condition
$$
	r(y,z,w) r(x,z,yw) r(x,y,w) s(x,y,z,w)
	= s'(x,y,z,w) r(x,y,zw) r(x,z,w) r(y,z,xw) \ ,
$$
that is, $s' = s \, dr$. 

Altogether we see that $H^4_b(S,B)$ describes b-bicategories of the form $\mathcal{B}(S,B,s)$ up to equivalences induced by b-functors whose underlying functor of bicategories is the identity. 

\section{Concluding remarks}\lb{concl}

Here we briefly list some questions which are left open in this paper and to which we would like to return in the future:
\begin{itemize}
\item
The relation between b-categories and braided monoidal categories should have a bicategorical version. It should be reasonably straightforward (if tedious) to extend the arguments of Section \ref{bmbc} to bicategories, i.e.\ to show that a braided monoidal bicategory has a b-bicategory structure. To extend  the constructions of Section \ref{ubcbmc} one first needs to define unital b-bicategory.
\item
It seems plausible that there should exist higher analogues of b-categories and b-bicategories, whose coherences are governed by higher dimensional permutohedra. If this is the case, their one object versions should correspond to higher dimensional versions of Zamolodchikov's tetrahedron equation.
\item
It is quite reasonable to expect that the comparison homomorphisms between abelian group cohomology and b-cohomology from Remark \ref{coho} should extend to all degrees.
\end{itemize}

\newcommand\arxiv[2]      {\href{http://arXiv.org/abs/#1}{#2}}
\newcommand\doi[2]        {\href{http://dx.doi.org/#1}{#2}}
\newcommand\httpurl[2]    {\href{http://#1}{#2}}

\end{document}